\documentclass[11pt]{article}

\usepackage{amssymb,amsmath,bm, bbm,graphicx}
\usepackage{algorithm}
\usepackage{algpseudocode}

\usepackage[colorlinks,
linkcolor=blue,
anchorcolor=blue,
citecolor=blue,
breaklinks=true,
unicode
]{hyperref}
\usepackage[nameinlink,capitalise]{cleveref}
\usepackage{natbib}

\usepackage{amsthm}
\newtheorem{theorem}{Theorem}
\crefname{theorem}{Theorem}{Theorems}
\newtheorem{assumption}{Assumption}
\crefname{assumption}{Assumption}{Assumptions}
\newtheorem{lemma}{Lemma}
\crefname{lemma}{Lemma}{Lemmas}

\crefname{example}{Example}{Examples}

\crefname{corollary}{Corollary}{corollaries}
\newtheorem{definition}{Definition}
\crefname{definition}{Definition}{Definitions}

\usepackage{enumitem}
\usepackage{appendix}
\usepackage{breakcites}

\newcommand*{\sign}{\mathrm{sign}}

\newcommand{\tr}{\textrm{trace}}

\newcommand{\Rom}[1]{\text{\uppercase\expandafter{\romannumeral #1\relax}}}

\newcommand\mc[1]{\mathcal{#1}}

\def\P{\mathbb{P}} 

\def\as{\lambda}
\def\rs{\lambda_r}

\providecommand{\tr}{\mathop\mathrm{tr}}
\newcommand{\opnorm}[1]{{\vert\kern-0.25ex\vert\kern-0.25ex\vert #1
    \vert\kern-0.25ex\vert\kern-0.25ex\vert}}
\DeclareMathOperator*{\esssup}{ess\,sup}

\def \RR{\mathbb{R}}
\def \cE{\mathcal{E}}
\def \BB{\mathbb{B}}
\def \PP{\mathbb{P}}
\def \EE{\mathbb{E}}

\newcommand{\scolor}[1]{{\color{black}#1}}
\newcommand{\bcolor}[1]{{\color{blue}#1}}

\newcommand{\scomment}[1]{\scolor{$\dagger$}\marginpar{\tiny\scolor{S:\ #1}}\hspace{-3pt}}

\usepackage{txfonts}

\usepackage{geometry}
\def\##1\#{\begin{align}#1\end{align}}
\def\$#1\${\begin{align*}#1\end{align*}}

\begin{document}

\title{ \LARGE Bernstein's inequalities for general Markov chains}


\author{Bai Jiang\thanks{Princeton University; Email: \url{baij@princeton.edu}.}  \and Qiang Sun\thanks{University of Toronto and MBZUAI; Email: \url{qsunstats@gmail.com}. His research is partially supported by the Natural Sciences and Engineering Research Council of Canada (Grant RGPIN-2018-06484) and by computing resources provided by the Digital Research Alliance of Canada.} \and Jianqing Fan\thanks{Princeton University; Email: \url{jqfan@princeton.edu}.  His research is supported by NSF grant DMS-2210833, DMS-2052926, and  ONR N00014-22-1-2340.}
}


\date{}

\maketitle


\begin{abstract}

We establish Bernstein's inequalities for functions of general (general-state-space and possibly non-reversible) Markov chains. These inequalities achieve sharp variance proxies and encompass the classical Bernstein inequality for independent random variables as special cases. The key analysis lies in bounding the operator norm of a perturbed Markov transition kernel by the exponential of sum of two convex functions. One coincides with what delivers the classical Bernstein inequality, and the other reflects the influence of the Markov dependence. A convex analysis on these two functions then derives our Bernstein inequalities. As applications, we apply our Bernstein inequalities to the Markov chain Monte Carlo integral estimation problem and the robust mean estimation problem with Markov-dependent samples, and achieve tight deviation bounds that previous inequalities can not. 
\end{abstract}

\tableofcontents

\section{Introduction} \label{section: introduction}

Concentration inequalities are tail probability bounds regarding how the average of random variables deviates from its expectation. They have found enormous applications in statistics, machine learning, information theory and other related fields. Many concentration inequalities are derived by the Cram\'{e}r-Chernoff method \citep{cramer1938mathematical,chernoff1952measure,boucheron2013concentration}. Let $Z_1, \dots, Z_n$ be $n$ random variables with mean zero. Suppose the moment generating function (mgf) of $\sum_{i=1}^n Z_i$ is bounded as
\begin{equation} \label{eqn1}
\mathbb{E}\left[e^{t\sum_{i=1}^n Z_i}\right] \le e^{ng(t)}, ~\forall~t \in \mathbb{R},
\end{equation}
with some convex function $g(t)$. Let $g^*(\epsilon) \coloneqq \sup \{t\epsilon - g(t): ~t \in \mathbb{R}\}$ for $\epsilon > 0$ be the Fenchel conjugate of $g(t)$, then the Cram\'{e}r-Chernoff method derives that
\begin{equation} \label{eqn2}
\mathbb{P}\left(\frac{1}{n}\sum_{i=1}^n Z_i \ge \epsilon\right) \le e^{- n g^*(\epsilon)},~\forall~\epsilon > 0.
\end{equation}
According to the conjugation and second-order properties of convex functions \citep{gorni1991conjugation}, if there exists a constant $V > 0$ such that
$$g(t) = \frac{Vt^2}{2} + o(t^2) ~~~\text{as $t \to 0$}, $$
then
$$g^*(\epsilon) = \frac{\epsilon^2}{2V} + o(\epsilon^2) ~~~ \text{as $\epsilon \to 0$}.$$
That means, a convex function $g(t)$ in the mgf bound \eqref{eqn1} with a sharp constant $V = \lim_{t \to 0} \frac{2g(t)}{t^2}$ leads to a tight concentration in the inequality \eqref{eqn2}.

We refer to this constant $V$ as the variance proxy of the concentration inequality, as it not only characterizes the shape of the probability tail in \eqref{eqn2} but also bounds the variance of $\bar{Z}_n := \frac{\sum_{i=1}^n Z_i}{\sqrt{n}}$. Indeed, bringing
\begin{align*}
\mathbb{E}\left[e^{t\sum_{i=1}^n Z_i}\right]
&= 1 + \mathbb{E}\left[\sum_{i=1}^n Z_i\right] \cdot t + \mathbb{E}\left[\left(\sum_{i=1}^n Z_i\right)^2 \right] \cdot \frac{t^2}{2} + \dots = 1 + \text{Var}\left[\sum_{i=1}^n Z_i \right] \cdot \frac{t^2}{2} + o(t^2),\\
e^{n g(t)}
&= 1 + n g(t) + \frac{(ng(t))^2}{2} + \dots = 1 + nV \cdot \frac{t^2}{2} + o(t^2)
\end{align*}
into both sides of \eqref{eqn1} yields $\text{Var}[\sqrt{n} \bar{Z}_n] \cdot t^2/2 \le nV \cdot t^2/2 + o(t^2)$. Dividing both sides by $nt^2/2$ and taking $t \to 0$ yield $\text{Var}[\bar{Z}_n] \le V$.

Hoeffding's, Bennett's and Bernstein's inequalities can be viewed as applications of the Cram\'{e}r-Chernoff method to independent random variables under boundedness and/or variance assumptions \citep{bernstein1946theory,bennett1962probability,hoeffding1963probability}. Suppose $Z_1,Z_2,\dots,Z_n$ are independent random variables with $\mathbb{E} Z_i = 0$ and $|Z_i| \le c$. Hoeffding's inequality holds in the form of \eqref{eqn1} and \eqref{eqn2} with
$$g(t) = \frac{c^2 t^2}{2}, ~~~ g^*(\epsilon) = \frac{\epsilon^2}{2c^2}, ~~~V = c^2.$$
Bennett's inequality takes variances of random variables into account and holds with
\begin{equation} \label{eqn:Bennett}
g(t) = \frac{\sigma^2}{c^2}(e^{tc} - 1 - tc), ~~~ g^*(\epsilon) = \frac{\sigma^2}{c^2}h\left(\frac{c\epsilon}{\sigma^2}\right), ~~~ V = \sigma^2,
\end{equation}
where $\sigma^2 = \frac{\sum_{i=1}^n \text{Var}(Z_i)}{n}$ and $h: u \in [0,\infty) \mapsto (1+u)\log(1+u) - u$. Lower bounding $h(u) \ge \frac{u^2}{2(1+u/3)}$ in Bennett's inequality further derives Bernstein's inequality
\begin{equation} \label{eqn:Bernstein}
\mathbb{P}\left(\frac{1}{n}\sum_{i=1}^n Z_i \ge \epsilon\right) \le \exp\left(-\frac{n\epsilon^2}{2(\sigma^2 + c\epsilon / 3)}\right).
\end{equation}
While Hoeffding's inequality uses the bound of random variables $|Z_i| \le c$ only, Bennett's and Bernstein's inequalities incorporate variances of $Z_i$ and thus are usually sharper when $\sigma$ is available and $\sigma \ll c$.

\scolor{Significant advances of concentration inequalities have been made, particularly with the work of \cite{talagrand1996new,ledoux2001concentration}. Their work laid the foundation for exploring concentration phenomena in more complex settings of various stochastic processes. While early research often focused on independent data, many real-world applications require the modeling of dependencies, providing a more realistic framework. Notable models of dependence include exchangeability and Markov chains. Exchangeability, particularly useful in the analysis of large networks, relies on de Finetti’s theorem \cite{definetti1937prevision} and the Aldous-Hoover-Kallenberg approach \citep{aldous1981representations,hoover1979relations,kallenberg2005probabilistic}. Recent studies, such as those in conformal prediction, have leveraged this modeling \citep{vovk2005algorithmic}. Markov chains, on the other hand, are frequently used in applications such as numerical integration via Monte Carlo methods \citep{gilks1995markov}, Hidden Markov Models \citep{baum1970maximization,rabiner1989introduction,bishop2006pattern}, and optimization algorithms like stochastic gradient descent \citep{bottou2010large,goodfellow2016deep} and simulated annealing \citep{kirkpatrick1983optimization}. Their exponential convergence to a stationary distribution, under spectral gap assumptions, is particularly valuable for deriving concentration inequalities. This has led to a growing body of work focused on understanding concentration in the context of Markov chains.}


Suppose that $Z_i = f_i(X_i)$, $1 \le i \le n$ are functions of a Markov chain $X_i$, $1 \le i \le n$ with a spectral gap. \cite{leon2004optimal} proved a Hoeffding-type inequality with a sharp variance proxy for a finite-state-space, reversible Markov chain and a time-independent function $f_i = f$. \cite{miasojedow2014hoeffding} obtained a variant of L\'eon and Perron's inequality for a general (general-state-space, possibly non-reversible) Markov chain. \citet*{fan2021hoeffding} improved upon them and established a Hoeffding-type inequality with the optimal variance proxy for general Markov chains and time-dependent functions $f_i$'s. \cite{lezaud1998chernoff} incorporated variances of $Z_i = f(X_i)$ into a Chernoff-type inequality for a finite-state-space, reversible Markov chain $\{X_i\}_{i \ge 1}$ and a time-independent function $f_i = f$. His analysis recently inspired several Bernstein-type inequalities in \citep{paulin2015concentration}.

A few other concentration inequalities for Markov chains are developed upon the minorization and drift conditions \citep{glynn2002hoeffding, douc2011consistency, adamczak2015exponential}, or information-theoretical ideas \citep{kontoyiannis2005relative, kontoyiannis2006exponential} rather than the spectral analysis. These inequalities therefore do not require the spectral gap assumption. But they usually have less explicit and often sub-optimal constants or more complicated expressions. A less related line of research includes Marton's work on the concentration of measure phenomenon for contracting Markov chains \citep{marton1996bounding, marton1996measure, marton1998measure, marton2003measure, marton2004measure} and further progresses under various dependence structures \citep{samson2000concentration, chazottes2007concentration, kontorovich2008concentration, redig2009concentration}.

This paper establishes two Bernstein-type inequalities for general Markov chains under spectral gap assumptions. One regards the time-dependent function case in which $f_i$, $i=1,2,\dots,n$ are different. This inequality encompasses the classical Bernstein's inequality for independent random variables. The other regards the time-independent function case in which $f_i = f$, $i=1,2,\dots,n$ are identical. This inequality encompasses the classical Bernstein's inequality for independently and identically distributed (i.i.d.) random variables. These two Bernstein-type inequalities achieve optimal variance proxies and tightest concentration, as their exponential rates are reached asymptotically by a class of Markov chains. 

The generality and sharpness of our inequalities stem from several new elements built upon Lezaud's method \citep{lezaud1998chernoff}. To outline the technical roadmap of this paper, we first revisit Lezaud's method. Suppose that $\{X_i\}_{i \ge 1}$ is a reversible, stationary Markov chain with invariant distribution $\pi$ on a finite state space $\mc{X}$, and that $f$ is a bounded and zero-mean function. Let $P$ denote the transition probability matrix of the Markov chain and $E^h$ denote the diagonal matrix with elements $e^{h(x)}$. Let $\opnorm{\cdot}_\pi$ denotes the operator norm in the real Hilbert space $\mc{L}_2 = \{h: \sum_{x \in \mc{X}} h^2(x)\pi(x) < \infty\}$. Lezaud's method consists of the following four steps.
\begin{enumerate}[label=(\alph*)]
\item For any $t > 0$, $\mathbb{E}_\pi\left[e^{t\sum_{i=1}^n f(X_i)}\right] \le \opnorm{E^{tf/2}PE^{tf/2}}_\pi^n$;
\item $\opnorm{E^{tf/2}PE^{tf/2}}_\pi$ is equal to the largest eigenvalue of $E^{tf/2}PE^{tf/2}$, which coincides with that of $PE^{tf}$;
\item By Kato's perturbation theory \citep[Chapter 2]{kato2013perturbation}, the largest eigenvalue of perturbed Markov operator $PE^{tf}$ is expanded as a series in $t$ and then upper bounded;
\item The Cram\'{e}r-Chernoff method derives final inequalities.
\end{enumerate}

Our first goal is to extend Lezaud's method for general-state-space Markov chains. The primary challenge is the inapplicability of step (b) to the general setting. In the general setting, $P$ is the Markov operator induced by the transition kernel, and $E^h$ is the multiplication order induced by function $e^{h(x)}$. The perturbed Markov operator $E^{tf/2}PE^{tf/2}$ may be infinite-dimensional and its operator norm is not necessarily equal to its largest eigenvalue. To address this challenge, we use techniques of L\'eon-Perron operator and function discretization, which were originally developed along the research line of Hoeffding-type inequalities for Markov chains \citep*{leon2004optimal,miasojedow2014hoeffding,fan2021hoeffding}. We take L\'eon-Perron operator $\widehat{P}$ as the extremal counterpart of the original Markov operator $P$ and replace step (a) of Lezaud's method with
\begin{align*}
\mathbb{E}_\pi\left[e^{t\sum_{i=1}^n f(X_i)}\right] \le \opnorm{E^{tf/2}\widehat{P} E^{tf/2}}_\pi^n.
\end{align*}
Next, by discretizing $f$ as $f_k$, a sequence of functions uniformly convergeing to $f$, we show that
$$\opnorm{E^{tf/2}\widehat{P} E^{tf/2}}_\pi = \lim_{k \to \infty} \opnorm{E^{tf_k/2}\widehat{P} E^{tf_k/2}}_\pi.$$
$E^{tf_k/2}\widehat{P}E^{tf_k/2}$ is self-adjoint and behaves like an operator on a finite state space.


We henceforth conduct spectral analyses of step (c) for $E^{tf_k/2}\widehat{P}E^{tf_k/2}$ in a finite state space. In this regard, we derive tighter estimates for series coefficients than Lezaud did in his original work. These tighter estimates lead to a more precise characterization of the mgf bound
$$\mathbb{E}_\pi\left[e^{t\sum_{i=1}^n f(X_i)}\right] \le e^{ng(t)},~~~g(t) = \underbrace{\frac{\sigma^2}{c^2} (e^{tc} - 1 - tc)}_{g_1(t)} + \underbrace{\frac{\sigma^2 \as t^2}{1- \as - 5ct}}_{g_2(t)},$$
where $c$ is the bound of $f$, $\sigma^2$ is the variance of $f$, $1-\lambda$ is the spectral gap of the Markov chain. In this bound, both $g_1(t)$ and $g_2(t)$ are convex functions of $t$. $g_1(t)$ coincides with the convex function that derives the classical Bernstein's inequality \eqref{eqn:Bennett} for independent random variables, and $g_2(t)$ goes to $0$ as the Markov dependence degenerates to the independence, i.e., $\lambda \to 0$. This preciseness eventually contributes to the exact reduction of our Bernstein inequalities to their classical counterparts under independence.

Another challenge arises when applying step (d) of Lezaud method. Due to the unavailability of an explicit expression of the convex conjugate of $g(t) = g_1(t)+g_2(t)$, we cannot easily convert the mgf bound to the deviation probability bound by using the Cram\'{e}r-Chernoff method. To address this challenge, we use the Moreau-Rockafellar formula \citep[Theorem 16.4]{rockafellar2015convex} in the convex analysis theory and lower bound the convex conjugate. This derivation loses no sharpness in terms of the variance proxy.



We provide applications of our newly derived Bernstein-type inequalities to \scolor{an MCMC integral estimation problem and a robust mean estimation problem. In the MCMC integral estimation problem, the non-asymptotic confidence interval by our Bernstein-type inequalities is at least $30\%$ tighter than previous state-of-the-art \citep{paulin2015concentration}.} For the robust mean estimation problem, the Hoeffding-type inequalities by \citet*{fan2021hoeffding} can only achieve a non-diminishing upper bound. This is because Bernstein-type inequalities incorporate variances of the random variables and thus are generally sharper, especially when the random variables  occasionally take on large values but have relatively small variances.


The rest of this paper proceeds as follows. \cref{sec2} presents our Bernstein-type inequalities and possible extensions to nonstationary Markov chains and unbounded random variables. \cref{sec3} summarizes preliminary results about the spectral analysis. \cref{sec4} collects technical proofs of main theorems. \cref{sec5} presents two applications of our newly-derived Bernstein-type inequalities. \cref{sec6} concludes the paper with a discussion.

\section{Main results}\label{sec2}

\subsection{Bernstein's inequalities for general Markov chains}

We introduce some notation before presenting main theorems. Denote by $\pi$ the stationary distribution of a Markov chain. Let $\pi(h) = \int h(x)\pi(dx)$ be the integral of real-valued function $h$ with respect to $\pi$. Denote by $\mc{L}_2 \coloneqq \{h: \pi(h^2) < \infty\}$ the real Hilbert space of $\pi$-square-integrable functions endowed with the inner product $\langle h_1, h_2\rangle_\pi = \pi(h_1h_2)$. Let $\mc{L}_2^0$ denote the subspace of $\mc{L}_2$ consisting of mean-zero functions. The transition kernel $P(x,y)$ of the Markov chain induces a Markov (integral) operator acting in $\mc{L}_2$, i.e., $Ph(x) = \int P(x,y)h(y)dy$.

\begin{definition}[Absolute spectral gap] \label{definition: absolute spectral gap}
A $\pi$-invariant Markov operator $P$ has non-zero absolute spectral gap $1-\as(P)$ if
$$\as(P) = \sup\left\{ \Vert Ph \Vert_\pi: \Vert h \Vert_\pi = 1, h \in \mc{L}_2^0\right\} < 1.$$
\end{definition}

\begin{definition}[Right spectral gap] \label{definition: right spectral gap}
A $\pi$-invariant Markov operator $P$ has non-zero right spectral gap $1-\rs(P)$ if
$$\rs(P) = \sup\left\{ \langle h, Ph \rangle_\pi: \Vert h \Vert_\pi = 1, h \in \mc{L}_2^0\right\} < 1.$$
\end{definition}

The next two theorems present the main results of this paper.

\begin{theorem} \label{thm1}
Suppose $\{X_i\}_{i \ge 1}$ is a stationary Markov chain with invariant distribution $\pi$ and absolute spectral gap $1-\as >0$, and $f_i$'s are functions with $|f_i| \le c$ and $\pi(f_i) = 0$. Let $\sigma^2 = \sum_{i=1}^n \pi(f_i^2) / n$. Then, for any $0 \le t < (1-\as)/5c$,
\begin{equation} \label{eqn:thm11}
\mathbb{E}_\pi\left[e^{t\sum_{i=1}^n f_i(X_i)}\right] \le e^{ng(t;\as, c, \sigma^2)}
\end{equation}
where
$$g(t; \as, c, \sigma^2) = \frac{\sigma^2}{c^2} (e^{tc} - 1 - tc) + \frac{\sigma^2 \lambda t^2}{1- \as - 5ct}.$$
Moreover, for any $\epsilon > 0$,
\begin{equation} \label{eqn:thm12}
\mathbb{P}_\pi\left(\frac{1}{n}\sum_{i=1}^n f_i(X_i) \ge \epsilon\right) \le e^{-ng^*(\epsilon; \as, c, \sigma^2)},
\end{equation}
where $g^*(\epsilon; \as, c, \sigma^2)$ is convex conjugate of $g(t; \as, c, \sigma^2)$ and
$$g^*(\epsilon; \as, c, \sigma^2) \ge \frac{n\epsilon^2/2}{\alpha_1(\as) \cdot \sigma^2 + \alpha_2(\as) \cdot c\epsilon}$$
with
$$\alpha_1(\as) = \frac{1+\as}{1-\as},~~~\alpha_2(\as) = \begin{cases}
1/3 &\mbox{if}~\as = 0,\\
5/(1-\as) &\mbox{if}~\as \in (0,1).
\end{cases}$$
\end{theorem}

\begin{theorem} \label{thm2}
Suppose $\{X_i\}_{i \ge 1}$ is a stationary, reversible Markov chain with invariant distribution $\pi$ and right spectral gap $1-\rs >0$, and $f$ is a function with $|f| \le c$ and $\pi(f) = 0$. Let $\sigma^2 = \pi(f^2)$. Then, for any $0 \le t < (1-\max\{\rs,0\})/5c$,
\begin{equation} \label{eqn:thm21}
\mathbb{E}_\pi\left[e^{t\sum_{i=1}^n f(X_i)}\right] \le e^{ng(t; \max\{\rs, 0\}, c, \sigma^2)}.
\end{equation}
Moreover, for any $\epsilon > 0$,
\begin{equation} \label{eqn:thm22}
\mathbb{P}_\pi\left(\frac{1}{n}\sum_{i=1}^n f(X_i) \ge \epsilon\right) \le e^{-ng^*(\epsilon; \max\{\rs, 0\}, c, \sigma^2)}.
\end{equation}
Both $g$ and $g^*$ are defined in \cref{thm1}.
\end{theorem}


These results compare well with the classical Bernstein's inequalities for independent random variables. In the mgf bound \eqref{eqn:thm11}, $g(t; \as, c, \sigma^2)$ is sum of two convex functions. The first convex function coincides with that in the mgf bound \eqref{eqn:Bennett} of the classical Bennett's and Bernstein's inequalities. The second convex function reflects the Markov dependence and decreases to $0$ as $\as$ decreases to zero. Consequently, \cref{thm1} encompasses the classical Bennett's and Bernstein's inequalities under independence as special cases. Indeed, independent random variables $Z_1,Z_2,\dots,Z_n$ can be seen as transformations of i.i.d. random variables $U_1,U_2,\dots,U_n \sim \text{Uniform}[0,1]$ via the inverse cumulative distribution functions $F_{Z_i}^{-1}$, i.e, $Z_i = f_i(U_i)$ with $f_i = F_{Z_i}^{-1}$. The i.i.d. random variables $\{U_i\}_{i \ge 1}$ form a stationary Markov chain with $\as = 0$. In this case, the mgf bound \eqref{eqn:thm11} in the Markov dependence scenario is simplified to be the mgf bound in the classical independence scenario. In a similar vein, \cref{thm2} encompasses the classical Bennett's and Bernstein's inequalities for i.i.d. random variables as special cases.

The variance proxy in \cref{thm1} is optimal in the sense that it is attainable by a class of Markov chains. Consider a reversible, stationary Markov chain $\{X_i\}_{i \ge 1}$ with transition kernel $P(x,B) = \as \mathbb{I}(x \in B) + (1-\as)\pi(B)$ with $\lambda \in (0,1)$ for state $x$ and any measurable subset $B$ of the state space. This chain has invariant distribution $\pi$ and absolute (and right) spectral gap $1-\as$. For function $f$ with $|f| \le c $, $\pi(f)=0$ and $\pi(f^2) = \sigma^2$, the Central Limit Theorem for reversible Markov chains \citep{kipnis1986central,geyer1992practical,rosenthal2003asymptotic} asserts that the asymptotic variance $\sigma^2_\text{as}(f)$ is given by
$$\sigma^2_\text{as}(f) \coloneqq \lim_{n \to \infty} \text{Var}[\bar{f}_n(X_1,\dots,X_n)] = \frac{1+\as}{1-\as} \cdot \sigma^2, ~~~\text{where}~ \bar{f}_n(X_1,\dots,X_n) = \frac{\sum_{i=1}^n f(X_i)}{\sqrt{n}}$$
Putting it together with the fact that any variance proxy $V \ge \text{Var}[\bar{f}_n]$ for any $n \ge 1$ yields
$$V \ge \sigma^2_\text{as}(f) = \lim_{n \to \infty} \text{Var}[\bar{f}_n] = \frac{1+\as}{1-\as} \cdot \sigma^2.$$
This lower bound for the variance proxy $V$ is attained by \eqref{eqn:thm12} in \cref{thm1}.

In a similar vein, the variance proxy of \eqref{eqn:thm22} in \cref{thm2} is optimal if $\rs \ge 0$, which is the case for Markov chains generated by popular MCMC algorithms such as the Metropolis-Hastings algorithm \citep{atchade2007geometric}. The assumption of $\rs > 0$ for reversible Markov chains is crucial for ensuring ergodicity, fast convergence, and a shorter mixing time. It facilitates theoretical analysis and provides a foundation for MCMC algorithms \citep{aldous1983random,diaconis1989geometric,robert1999monte,levin2009markov}. However, there do exist reversible Markov chains with $\rs < 0$, which typically indicates a oscillatory behavior. For example,
$$P = \begin{pmatrix}
0.5 & 0.4 & 0.1 \\
0.4 & 0.2 & 0.4 \\
0.1 & 0.4 & 0.5
\end{pmatrix},~~~\rs(P) = -0.2.$$

\cref{table:proxy1} summarizes some of the recent results on Hoeffding-type and Bernstein-type inequalities for Markov chains along with the classical Hoeffding's, Bennett's and Bernstein's inequalities for the time-dependent function case. \cref{table:proxy2} summarizes results for the time-independent function case. Notation in both tables are the same as in \cref{thm1,thm2}. As shown in \cref{table:proxy1}, \cref{thm1} in this paper improves over the previous Bernstein-type inequality of \citet[Equation 3.22 in Theorem 3.9]{paulin2015concentration} in terms of the sharpness of variance proxy. Specifically, $\frac{1+\as}{1-\as} < \frac{4}{1-\as^2}$ for any $\as \in [0,1).$

\begin{table}[t!]
\caption{Concentration inequalities for time-dependent functions $f_i$ of Markov chains.}
\label{table:proxy1}
\begin{tabular}{llll}
\hline
type & reference & condition & variance proxy\\
\hline
\\[-0.8em]
Hoeffding & \citealp{hoeffding1963probability} & independent & $c^2$ \\
\\[-0.8em]
Hoeffding & \citealp{fan2021hoeffding} & general-state-space & $\frac{1+\as}{1-\as} \cdot c^2$\\
\\[-0.8em]
Bennett & \citealp{bennett1962probability} & independent & $\sigma^2$ \\
\\[-0.8em]
Bernstein & \citealp{bernstein1946theory} & independent & $\sigma^2$ \\
\\[-0.8em]
Bernstein & \citealp{paulin2015concentration}, (3.22) & general-state-space, reversible\footnote{Who am I} & $\frac{4}{1-\as^2} \cdot \sigma^2$ \\
\\[-0.8em]
Bernstein & \cref{thm1} & general-state-space & $\frac{1+\as}{1-\as} \cdot \sigma^2$\\
\\[-0.8em]
\hline
\end{tabular}
\begin{flushleft}
\footnotesize \footnotemark[1] \citet{paulin2015concentration} defines spectral gap for non-reversible Markov chains in a different way to ours. His Bernstein-type inequalities for non-reversible Markov chains involve a quantity called ``pseudo spectral gap''. The relation between his ``pseudo spectral gap'' and our ``absolute spectral gap'' is unclear, and thus this table does not include his Bernstein-type inequalities for non-reversible Markov chains.
\end{flushleft}
\end{table}

\begin{table}
\caption{Concentration inequalities for time-independent function $f$ of Markov chains.}
\label{table:proxy2}
\begin{tabular}{llll}
\hline
type & reference & condition & variance proxy\\
\hline
\\[-0.8em]
Hoeffding & \citealp{hoeffding1963probability} & independent & $c^2$ \\
\\[-0.8em]
Hoeffding & \citealp{miasojedow2014hoeffding} & general-state-space & $\frac{1+\as}{1-\as} \cdot c^2$ \\
\\[-0.8em]
Hoeffding & \citealp{fan2021hoeffding} & general-state-space & $\frac{1+\as}{1-\as} \cdot c^2$ \\
\\[-0.8em]
Hoeffding & \citealp{leon2004optimal} & finite-state-space, reversible & $\frac{1+\max\{\rs,0\}}{1-\max\{\rs,0\}} \cdot c^2$ \\
\\[-0.8em]
Hoeffding & \citealp{fan2021hoeffding} & general-state-space, reversible & $\frac{1+\max\{\rs,0\}}{1-\max\{\rs,0\}} \cdot c^2$ \\
\\[-0.8em]
Bennett & \citealp{bennett1962probability} & independent & $\sigma^2$ \\
\\[-0.8em]
Bernstein & \citealp{bernstein1946theory} & independent & $\sigma^2$ \\
\\[-0.8em]
Chernoff & \citealp{lezaud1998chernoff}, (13) & general-state-space & $\frac{4}{1-\as} \cdot \sigma^2$ \\
\\[-0.8em]
Bernstein & \cref{thm1} & general-state-space & $\frac{1+\as}{1-\as} \cdot \sigma^2$\\
\\[-0.8em]
Chernoff & \citealp{lezaud1998chernoff}, (1) & finite-state-space, reversible & $\frac{2}{1-\rs} \cdot \sigma^2$ \\
\\[-0.8em]
Bernstein & \citealp{paulin2015concentration}, (3.20)\footnote{Who am I} & general-state-space, reversible & $\left(\frac{1+\rs}{1-\rs} + 0.8\right) \cdot \sigma^2$\\
\\[-0.8em]
Bernstein & \citealp{paulin2015concentration}, (3.21) & general-state-space, reversible & $\frac{2}{1-\rs} \cdot \sigma^2$ \\
\\[-0.8em]
Bernstein & \cref{thm2} & general-state-space, reversible & $\frac{1+\max\{\rs,0\}}{1-\max\{\rs,0\}} \cdot \sigma^2$\\
\\[-0.8em]
\hline
\end{tabular}
\begin{flushleft}
\footnotesize \footnotemark[2] (3.20) of \cite{paulin2015concentration} gives a variance proxy  $\sigma^2_\text{as}(f) + 0.8\sigma^2$. Since $\sigma^2_\text{as}(f)$ is usually unknown in practice and is equal to $\frac{1+\rs}{1-\rs} \cdot \sigma^2$ in the worst case for reversible Markov chains \citep[Proposition 1]{rosenthal2003asymptotic}, we replace $\sigma^2_\text{as}(f)$ with $\frac{1+\rs}{1-\rs} \cdot \sigma^2$ for clearness.
\end{flushleft}
\end{table}

\subsection{Extendable and non-extendable}\label{sec:ext}

In this subsection, we extend our Bernstein-type inequalities to the nonstationary case and present that our Bernstein-type inequalities are not extendable to unbounded functions. 

We first extend \cref{thm1,thm2} to non-stationary Markov chains, provided that the initial distribution of the Markov chain  exhibits finite moments with respect to the invariant distribution. It is noteworthy that the trick of this extension is an elementary use of measure changing and H\"{o}lder's inequality, which should be applicable to any Bernstein-type inequality in a similar setup.

\begin{theorem} \label{thm3}
Suppose $\{X_i\}_{i \ge 1}$ is a Markov chain with initial distribution $\nu$, invariant distribution $\pi$ and absolute spectral gap $1-\as >0$, and $f_i$ are functions with $|f_i| \le c$ and $\pi(f_i) = 0$. Let $\sigma^2 = \sum_{i=1}^n \pi(f_i^2) / n$. Suppose $\nu$ is absolutely continuous with respect to $\pi$ and its derivative $\frac{d\nu}{d\pi}$ has a finite $p$-moment for some $p \in (1,\infty]$, i.e.,
$$\left\Vert \frac{d\nu}{d\pi} \right\Vert_{\pi,p} \coloneqq \begin{cases}
\left[\int \left|\frac{d\nu}{d\pi}(x)\right|^p \pi(dx)\right]^{1/p} & \mbox{if}~ 1 < p < \infty\\
\esssup \left|\frac{d\nu}{d\pi}(x)\right| & \mbox{if}~p = \infty\\
\end{cases} < \infty
$$
Let $q = p/(p-1) \in [1, \infty)$ then, for any $0 \le t < (1-\as)/5cq$,
\begin{equation*}
\mathbb{E}_{\nu}\left[e^{t\sum_{i=1}^{n} f(X_i)}\right] \le \left\Vert \frac{d\nu}{d\pi} \right\Vert_{\pi,p} e^{n g(qt) / q}
\end{equation*}
It follows that
\begin{equation*}
\mathbb{P}_\nu \left( \frac{1}{n}\sum_{i=1}^{n} f_i(X_i) \ge \epsilon \right) \le \left\Vert \frac{d\nu}{d\pi} \right\Vert_{\pi,p}  e^{-ng^*(\epsilon)/q},
\end{equation*}
where $g(t) = g(t; \as, c, \sigma^2)$ and $g^*(\epsilon) = g^*(\epsilon; \as, c, \sigma^2)$ are defined in \cref{thm1}.
If the Markov chain is reversible and $f_i = f$ for $i=1,\dots,n$, then the above two inequalities holds with $\max\{\rs, 0\}$ in place of $\as$.
\end{theorem}

The classical Bernstein's inequality can be extended for some unbounded (e.g., sub-Exponential) independent random variables. Unfortunately, the boundedness of random variables is necessary for our Bernstein-type inequalities to hold in the Markov chain setup. Below is a strengthened version of \citet[Theorem 4]{fan2021hoeffding}. It shows that convex function $g(t) = Vt^2/2 + o(t^2)$ is not available in an unbounded case.

\begin{theorem}\label{thm4}
Consider a stationary Markov chain $\{X_i\}_{i \ge 0}$ on the state space $\mathcal{X}=\mathbb{R}$ with invariant distribution $\pi \sim \mathcal{N}(0,1)$ and transition kernel
$$P(x,B) = \as \mathbb{I}(x \in B) + (1-\as)\pi(B), ~~~ \forall x \in \mathcal{X}, ~\forall~\text{measurable}~B \subseteq \mathcal{X},$$
for some $0 < \as < 1$. There exists no function $g(t) = Vt^2/2 + o(t^2)$ with finite variance proxy $V>0$ and $t_0 > 0$ such that 
$$\mathbb{E}_{\pi}{e^{t\sum_{i=1}^n X_i}} \le e^{ng(t)}, ~~~\forall n \ge 1, \forall t \in [0,t_0).$$
\end{theorem}

\section{Preliminaries}\label{sec3}
This section collects preliminary results needed for the technical proofs in this paper. Assume the state space of the Markov chain, denoted by $\mathcal{X}$, is equipped with a sigma-algebra $\mathcal{B}$, and ($\mathcal{X}, \mathcal{B})$ a standard Borel space\footnote{A measurable space $(\mathcal{X},\mathcal{B})$ is standard Borel if it is isomorphic to a subset of $\mathbb{R}$. See Definition 4.33 in \cite{breiman1992probability}.}. This assumption is commonly seen in the literature for measure-theoretical studies of Markov chains. It holds in most practical examples, in which $\mathcal{X}$ is a subset of a multi-dimensional real space and $\mathcal{B}$ is the Borel sigma-algebra over $\mathcal{X}$. Let $\{X_i\}_{i \ge 1}$ be a Markov chain on the state space $(\mathcal{X},\mathcal{B})$ with invariant measure $\pi$.

\subsection{Hilbert space and spectral gaps}
For a real-valued, $\mathcal{B}$-measurable function $h: \mathcal{X} \to \mathbb{R}$, let
$$\pi(h) \coloneqq \int_\mathcal{X} h(x)\pi(dx).$$
The set of all $\pi$-square-integrable functions
$$\mc{L}_2(\mathcal{X}, \mathcal{B}, \pi) \coloneqq \left\{ h: \pi(h^2) < \infty \right\}$$
is a real Hilbert space endowed with the inner product
$$\langle h_1, h_2 \rangle_\pi = \int_\mathcal{X} h_1(x)h_2(x)\pi(dx), ~~~\forall h_1,h_2 \in \mc{L}_2(\mathcal{X}, \mathcal{B}, \pi).$$
Define the norm of a function $h \in \mc{L}_2(\mathcal{X},\mathcal{B},\pi)$  as
$$\Vert h \Vert_\pi = \sqrt{\langle h, h \rangle_\pi},$$
which induces the norm of a linear operator $T$ on $\mc{L}_2(\mathcal{X},\mathcal{B},\pi)$ as
$$\opnorm{T}_\pi = \sup \{ \Vert Th \Vert_\pi: \Vert h \Vert_\pi = 1 \}.$$
Write $\mc{L}_2$ in place of $\mc{L}_2(\mathcal{X}, \mathcal{B}, \pi)$ for simplicity, whenever the probability space $(\mathcal{X}, \mathcal{B}, \pi)$ is clear in the context. Let $\mc{L}_2^0 = \{h \in \mc{L}_2:~ \pi(h) = 0\}$ denote the subspace of $\mc{L}_2$ consisting of mean-zero functions.

Each transition kernel $P(x,B)$ for $x \in \mathcal{X}$ and $B \in \mathcal{B}$, if invariant with respect to $\pi$, corresponds to an integral operator $h \mapsto \int_{\mathcal{X}} h(y)P(\cdot,dy)$ on $\mc{L}_2$. We abuse notation $P$ to denote this operator. Let $\mathit{1}: x \in \mathcal{X} \mapsto 1$ be identically one function, and let $\mathit{\Pi}$ be the projection operator onto $\mathit{1}$, i.e.,
$\mathit{\Pi}: ~h \mapsto \langle h, \mathit{1}\rangle_\pi \mathit{1} = \pi(h)\mathit{1}$. Then an equivalent definition of $\as(P)$ in \cref{definition: absolute spectral gap} is given by
$$\as(P) \coloneqq \opnorm{P - \mathit{\Pi}}_\pi < 1.$$
 
Denote by $S_0(P)$ the spectrum of operator $P$ acting in (complexified) $\mc{L}_2^0$. If $P$ is self-adjoint, or equivalently the Markov chain is reversible, then $S_0(P) \subseteq [-1, +1]$ are real-valued. In this case, an equivalent definition of $\rs(P)$ in \cref{definition: right spectral gap} is given by
$$\rs(P) = \sup\{s: s \in S_0(P)\}.$$

\subsection{Le\'on-Perron: convex majorization of Markov operator}
Note that every convex combination of Markov operators is still a Markov operator.

\begin{definition}[Le\'on-Perron operator]
A Markov operator $\widehat{P}_\lambda$ on $\mc{L}_2$ is said Le\'on-Perron if it is a convex combination of operators $I$ and $\mathit{\Pi}$ with some coefficient $\lambda \in [0,1)$, that is
$$\widehat{P}_\lambda = \lambda I + (1-\lambda)\mathit{\Pi}.$$
\end{definition}
This Le\'on-Perron operator associates with a transition kernel
$$\widehat{P}_\lambda(x,B) = \lambda \mathbb{I}(x \in B) + (1-\lambda) \pi(B), ~~~\forall x \in \mathcal{X}, \ \forall B \in \mathcal{B},$$
which characterizes a simple transition dynamics: at each step, the Markov chain either stays at the current state with probability $\lambda$ or jumps to a new state drawn from $\pi$ with probability $1-\lambda$. Evidently, any Le\'on-Perron operator $\widehat{P}_\lambda$ is invariant with $\pi$.


The Le\'on-Perron operator has played a central role as the convex majorization of the Markov operator in the argument of \citet*{leon2004optimal}. This paper uses a few properties of Le\'on-Perron operators for general Markov chains, which were developed along the line of research on Hoeffding-type inequalities for Markov chains \citep*{leon2004optimal,miasojedow2014hoeffding,fan2021hoeffding}. 

Let $E^h$ denote the multiplication operator of function $e^h: x \mapsto e^{h(x)}$. 
The following three lemmas are taken from \cite*{fan2021hoeffding}.
\begin{lemma}
\label{lemma:norm-product1}
Suppose $\{X_i\}_{i \ge 1}$ is a stationary Markov chain with invariant distribution $\pi$ and absolute spectral gap $1-\as>0$. For any bounded functions $f_i$ and any $t \in \mathbb{R}$,
$$\mathbb{E}_\pi\left[e^{t\sum_{i=1}^n f_i(X_i)}\right] \le \prod_{i=1}^n \opnorm{E^{tf_i/2}\widehat{P}_\as E^{tf_i/2}}_\pi.$$
\end{lemma}

\begin{lemma}
\label{lemma:norm-product2}
Suppose $\{X_i\}_{i \ge 1}$ is a stationary, reversible Markov chain with invariant distribution $\pi$ and right spectral gap $1-\rs>0$. For any bounded function $f$ and any $t \in \mathbb{R}$,
$$\mathbb{E}_\pi\left[e^{t\sum_{i=1}^n f(X_i)}\right] \le \opnorm{E^{tf/2}\widehat{P}_{\max\{\rs,0\}}E^{tf/2}}_\pi^n.$$
\end{lemma}

\begin{lemma}
\label{lemma:large deviation}
Suppose $\{\widehat{X}_i\}_{i \ge 1}$ is a stationary Markov chain driven by a Le\'on-Perron operator $\widehat{P}_\lambda = \lambda I + (1-\lambda)\mathit{\Pi}$ for some $\lambda \in [0,1)$. For any bounded function $f$ and any $t \in \mathbb{R}$,
$$\lim_{n \to \infty} \frac{1}{n} \log \mathbb{E}_\pi\left[e^{t\sum_{i=1}^n f(\widehat{X}_i)}\right] = \log \opnorm{E^{tf/2}\widehat{P}_\lambda E^{tf/2}}_\pi.$$
\end{lemma}

\cref{lemma:norm-product1} bounds the mgf of $\sum_{i=1}^n f_i(X_i)$ by the product of the operator norms of symmetrically perturbed Le\'{o}n-Perron operators $E^{tf_i/2}\widehat{P}_\as E^{tf_i/2}$. \cref{lemma:norm-product2} refines \cref{lemma:norm-product1} with $\max\{\rs,0\}$ in place of $\as$ in case that the Markov chain is reversible and $f_i = f$. \cref{lemma:large deviation} asserts that a Markov chain driven by $\widehat{P}_\lambda = \lambda I + (1-\lambda)\mathit{\Pi}$ is the extreme case of all Markov chains with absolute spectral gap $1-\lambda$ in the sense that the inequality in \cref{lemma:norm-product1} is asymptotically tight.

\subsection{Kato's perturbation theory}

Another cornerstone of our technique is Kato's analysis on the largest eigenvalue of a perturbed operator \citep{kato2013perturbation}. This was the main tool to prove Chernoff-type and Bernstein-type inequalities for Markov chains by \citet{lezaud1998chernoff} and \citet{paulin2015concentration}. Specifically, Kato's analysis expends the largest eigenvalue of the perturbed Markov operator $PE^{tf}$ as a series in $t$.

Consider a reversible Markov chain $\{X_i\}_{i \ge 1}$ with invariant distribution $\pi$, self-adjoint Markov operator $P$ and right spectral gap $1-\rs$. Recall that, for any function $f \in \mc{L}_2^0$, its asymptotic variance is given by
$$\sigma^2_\text{as}(f) = \lim_{n \to \infty} \text{Var}[\bar{f}_n(X_1,\dots,X_n)], ~~~ \bar{f}_n(X_1,\dots,X_n) = \frac{\sum_{i=1}^n f(X_i)}{\sqrt{n}}.$$
By \citet[Proposition 1]{rosenthal2003asymptotic}, for any $h \in \mc{L}_2^0$,
$$\sigma^2_\text{as}(f) \le \frac{1+\rs}{1-\rs} \cdot \Vert f\Vert_\pi^2.$$

Let $Z = (I - P + \mathit{\Pi})^{-1} -\mathit{\Pi}$ be the negative of the reduced resolvent of $P$ with respect to its largest eigenvalue $1$. \cref{lemma:resolvent} summarizes a few useful properties of $Z$, whose proof can be found in \citet[Proposition 1.5]{lezaud1998thesis} and \citet[Lemma 5.3]{paulin2015concentration}. \cref{lemma:lezaud}, borrowed from \citet{lezaud1998chernoff}, expends the largest eigenvalue of the perturbed Markov operator $PE^{tf}$ as a series in $t$ involving $Z$.

\begin{lemma} \label{lemma:resolvent}
Suppose a self-adjoint Markov operator $P$ has invariant distribution $\pi$ and right spectral gap $1-\rs$. The following properties of $Z = (I - P + \mathit{\Pi})^{-1} - \mathit{\Pi}$ hold.
\begin{enumerate}[label=(\alph*)]  
    \item $Z \mathit{\Pi} = \mathit{\Pi} Z = 0$.
    \item $ZP = Z - (I - \mathit{\Pi})$.
    \item $\langle h, Zh \rangle_\pi = (\sigma^2_\text{as}(h) + \Vert h\Vert_\pi^2) / 2,~~~\forall~h \in \mc{L}_2^0$.
    \item $\opnorm{Z}_\pi = (1-\rs)^{-1}$.
\end{enumerate}
\end{lemma}

\begin{lemma} \label{lemma:lezaud}
Consider a reversible, irreducible Markov chain on finite state space $\mathcal{X}$ with invariant distribution $\pi$, self-adjoint Markov operator (transition probability matrix) $P$ and right spectral gap $1-\rs > 0$ (in this case, $\rs$ is the second largest eigenvalue of $P$). Let $D$ be the diagonal matrix with elements $\{f(x): x \in \mathcal{X}\}$ and let $T^{(m)} = PD^m / m!$ for any $m \ge 0$ with $D^0 = I$ by convention. Then
$$PE^{tf} = P\left(\sum_{m=0}^\infty \frac{D^m}{m!} \cdot t^m\right) = \sum_{m=0}^\infty T^{(m)}\cdot t^m.$$
Let $t_0 =\left(2\opnorm{T^{(1)}}_\pi(1-\rs)^{-1} + c_0\right)^{-1}$ for some $c_0$ such that 
$$\opnorm{T^{(m)}}_\pi \le \opnorm{T^{(1)}}_\pi c_0^{m-1},~\forall m \ge 1.$$
Then, for any $t$ such that $|t| < t_0$, the largest eigenvalue of $PE^{tf}$, denoted by $\beta(t)$, admits an expansion 
$$\beta(t) = \sum_{m=0}^\infty \beta^{(m)}t^m.$$
In this expansion, $\beta^{(0)} = 1$ is the largest eigenvalue of $T^{(0)} = P$ and, for any $m \ge 1$,
\begin{equation*}
\beta^{(m)} = \sum_{p=1}^m \frac{-1}{p} \sum_{\begin{split}
v_1+\dots+v_p = m, ~v_i \ge 1\\
k_1+\dots+k_p = p-1,~k_j \ge 0\end{split}} \tr\left(T^{(v_1)}Z^{(k_1)}\dots T^{(v_p)}Z^{(k_p)}\right),
\end{equation*}
where $Z^{(0)} = - \mathit{\Pi}$ and $Z^{(j)} = Z^j$ for $j \ge 1$ are powers of $Z$.
\end{lemma}

\section{Proofs of Theorems}\label{sec4}

This section proves four lemmas (\cref{lemma:discretization,lemma:counterpart,lemma:finite-state-space,lemma:infimal convolution}) and then  main theorems.

\begin{lemma} \label{lemma:discretization}
Suppose $f_k$ is a sequence of functions uniformly converging to a bounded function $f$. For any bounded linear operator $T$ acting on $\mc{L}_2$ and any $t \in \mathbb{R}$,
$$\opnorm{E^{tf/2} T E^{tf/2}}_\pi = \lim_{k \to \infty} \opnorm{E^{tf_k/2} T E^{tf_k/2}}_\pi.$$
\end{lemma}
\begin{proof}[Proof of \cref{lemma:discretization}]
Let $\epsilon_k = \sup_x |f_k(x) - f(x)|$. Write
\begin{align*}
\opnorm{E^{tf/2} T E^{tf/2}}_\pi &= \opnorm{E^{t(f - f_k)/2} E^{tf_k/2} T E^{tf_k/2} E^{t(f - f_k)/2}}_\pi\\
&\le \opnorm{E^{t(f - f_k)/2}}_\pi \times \opnorm{E^{tf_k/2} T E^{tf_k/2}}_\pi \times \opnorm{E^{t(f - f_k)/2}}_\pi\\
&\le e^{t\epsilon_k}  \opnorm{E^{tf_k/2} T E^{tf_k/2}}_\pi.
\end{align*}
Similarly,
\begin{align*}
\opnorm{E^{tf_k/2} T E^{tf_k/2}}_\pi &\le e^{t\epsilon_k} \opnorm{E^{tf/2} T E^{tf/2}}_\pi.
\end{align*}
Putting the last two displays together yields
$$e^{-t\epsilon_k} \opnorm{E^{tf/2} T E^{tf/2}}_\pi  \le \opnorm{E^{tf_k/2} T E^{tf_k/2}}_\pi \le e^{+t\epsilon_k}  \opnorm{E^{tf/2} T E^{tf/2}}_\pi.$$
Taking $k \to \infty$ completes the proof.
\end{proof}

\begin{lemma} \label{lemma:counterpart}
Let $\widehat{P}_\lambda = \lambda I + (1-\lambda) \mathit{\Pi}$ be a Le\'on-Perron operator with $\lambda \in [0,1)$ on a general state space $\mathcal{X}$. Let $f$ be a function on $\mathcal{X}$ taking finitely many possible values. On the finite state space
$$\mathcal{Y} = \{ y \in f(\mathcal{X}): \pi(\{x:f(x)=y\})>0\},$$
define a transition probability matrix $\widehat{Q}_\lambda = \lambda I + (1-\lambda) \bm{1} \mu^\mathrm{T}$, with vector $\mu$ consisting of elements $\pi(\{x: f(x) = y\})$ for each $y \in \mathcal{Y}$. Let $E^{t\mathcal{Y}}$ denote the diagonal matrix with elements $\{e^{ty}:~y \in \mathcal{Y}\}$. Then
$$\opnorm{E^{tf/2}\widehat{P}_\lambda E^{tf/2}}_\pi = \opnorm{E^{t\mathcal{Y}/2}\widehat{Q}_\lambda E^{t\mathcal{Y}/2}}_\mu.$$
\end{lemma}

\begin{proof}[Proof of \cref{lemma:counterpart}]
Let $\{B_i\}_{i \ge 1}$ be a sequence of i.i.d. Bernoulli random variables with success probability $\lambda$, and let $\{W_i\}_{i \ge 1}$ be a sequence of i.i.d. random variables following $\pi$, respectively. Evidently,
\begin{align*}
\widehat{X}_1 = W_1, ~~~ &\widehat{X}_i = B_i \widehat{X}_{i-1} + (1-B_i) W_i, ~~~\forall i \ge 2;\\
\widehat{Y}_1 = f(W_1), ~~~& \widehat{Y}_i = B_i \widehat{Y}_{i-1} + (1-B_i) f(W_i), ~~~\forall i \ge 2.
\end{align*}
are stationary Markov chains with invariant distributions $\pi, \mu$ and transition probabilities $\widehat{P}_\lambda, \widehat{Q}_\lambda$ on state spaces $\mathcal{X}, \mathcal{Y}$, respectively; and $\widehat{Y}_i = f(\widehat{X}_i)$ for $i \ge 1$. Putting them with \cref{lemma:large deviation} together yields
\begin{align*}
\log \opnorm{E^{tf/2}\widehat{P}_\lambda E^{tf/2}}_\pi
&= \lim_{n \to \infty} \frac{1}{n} \log \mathbb{E}_\pi \left[\exp \left(t\sum_{i=1}^n f(\widehat{X}_i)\right)\right]\\
&= \lim_{n \to \infty} \frac{1}{n} \log \mathbb{E}_\mu \left[\exp \left(t\sum_{i=1}^n \widehat{Y}_i\right)\right]\\
&= \log \opnorm{E^{t\mathcal{Y}/2}\widehat{Q}_\lambda E^{t\mathcal{Y}/2}}_\mu.
\end{align*}
\end{proof}

\begin{lemma} \label{lemma:finite-state-space}
Consider a reversible, irreducible Markov chain on finite state space $\mathcal{X}$ with invariant distribution $\pi$, self-adjoint Markov operator $P$, absolute spectral gap $1 - \as > 0$ and right spectral gap $1-\rs > 0$. For any function $f$ such that $|f| \le c$, $\pi(f) = 0$, and $\pi(f^2) = \sigma^2$ and any $0 \le t < (1-\rs)/5c$,
$$\opnorm{E^{tf/2}PE^{tf/2}}_\pi \le \exp{\left(\frac{\sigma^2}{c^2} (e^{tc} - 1 - tc) + \frac{\sigma^2 \as t^2}{1- \rs - 5ct}\right)}.$$
\end{lemma}

\begin{proof}[Proof of \cref{lemma:finite-state-space}]
Elements in the matrix $E^{tf/2}PE^{tf/2}$ are non-negative. It follows from the Perron-Frobenius theorem that the operator norm of $E^{tf/2}PE^{tf/2}$ coincides with its largest eigenvalue. $E^{tf/2}PE^{tf/2}$ is similar to $PE^{tf}$, so they share the same eigenvalues. Thus $\opnorm{E^{tf/2}PE^{tf/2}}_\pi$ is equal to the largest eigenvalue of $PE^{tf}$. Denote by $\beta(t)$ this eigenvalue. Recall from \cref{lemma:resolvent,lemma:lezaud} that
\begin{itemize}
    \item $Z = (I - P + \mathit{\Pi})^{-1} - \mathit{\Pi}$,
    \item $Z^{(0)} = -\Pi$, $Z^{(j)} = Z^j$ for $j \ge 1$.
    \item $D$ is the diagonal matrix with elements $\{f(x): x \in \mc{X}\}$,
    \item $T^{(0)} = P$, $T^{(m)} = PD^m/m!$ for $j \ge 1$.
\end{itemize}
\cref{lemma:lezaud} asserts that, for some constant $t_0 > 0$,
$$\beta(t) = \sum_{m=0}^\infty \beta^{(m)}t^m, ~~~\forall 0 \le t < t_0,$$
where $\beta^{(0)} = 1$ and
\begin{equation*}
\beta^{(m)} = \sum_{p=1}^m \frac{1}{p} \sum_{\begin{split}
v_1+\dots+v_p = m, ~v_i \ge 1\\
k_1+\dots+k_p = p-1,~k_j \ge 0\end{split}} -\frac{\tr\left(PD^{v_1}Z^{(k_1)}\dots PD^{v_p}Z^{(k_p)}\right)}{v_1!\dots v_p!}, ~~~\forall~m \ge 1.
\end{equation*}
The choice of $t_0 = \frac{1-\rs}{3-\rs} \times c^{-1}$ meets the requirement of \cref{lemma:lezaud}, as
\begin{align*}
\opnorm{T^{(1)}}_\pi &= \opnorm{PD}_\pi \le \opnorm{P}_\pi \opnorm{D}_\pi \le c,\\
\opnorm{T^{(m)}}_\pi &= \frac{1}{m!} \opnorm{PD^m}_\pi \le \opnorm{PD}_\pi \opnorm{D}_\pi^{m-1} \le \opnorm{T^{(1)}}_\pi c^{m-1}.
\end{align*}
Proceed to compute coefficients $\beta^{(m)}$ for $m \ge 1$. It is straightforward that
$$\beta^{(1)} = -\tr\left(PD^1Z^{(0)}\right) = \tr\left(PD\mathit{\Pi}\right) = \tr\left(D\mathit{\Pi}P\right) = \tr\left(D\mathit{\Pi}\right) = \pi(f) = 0.$$
By parts (b) and (c) of \cref{lemma:resolvent},
\begin{align*}
\beta^{(2)}
&= -\frac{1}{1}\frac{\tr\left(PD^2Z^{(0)}\right)}{2!} -\frac{1}{2}\frac{\tr\left(PD^1Z^{(0)}PD^1Z^{(1)}\right) + \tr\left(PD^1Z^{(1)}PD^1Z^{(0)}\right)}{1!1!}\\
&= \frac{\tr\left(PD^2\mathit{\Pi}\right)}{2} + \frac{\tr\left(PD\mathit{\Pi}PDZ\right) + \tr\left(PDZPD\mathit{\Pi}\right)}{2}\\
&= \frac{\tr\left(D^2\mathit{\Pi}P\right)}{2} + \tr\left(ZPD\mathit{\Pi}PD\right)\\
&= \frac{\tr\left(D^2\mathit{\Pi}\right)}{2} + \tr\left(ZPD\mathit{\Pi}D\right)\\
&= \frac{\Vert f \Vert_\pi^2}{2} + \langle ZPf, f\rangle_\pi\\
&= \frac{\sigma^2_\text{as}(f)}{2}.
\end{align*}
For $m \ge 3$,  $\beta^{(m)}$ is the sum of $m$ terms $\{\beta^{(m)}_p: ~p=1,\dots,m\}$.
$$\beta^{(m)}_p = \frac{1}{p} \sum_{\begin{split}
v_1+\dots+v_p = m, ~v_j \ge 1\\
k_1+\dots+k_p = p-1,~k_j \ge 0\end{split}} -\frac{\tr\left(PD^{v_1}Z^{(k_1)}\dots PD^{v_p}Z^{(k_p)}\right)}{v_1!\dots v_p!}.$$
For $p=1$,
$$\beta^{(m)}_1 = -\frac{\tr\left(PD^{m}Z^{(0)}\right)}{m!} = \frac{\pi(f^m)}{m!}.$$
For $p=2,\dots,m$, consider each term in the summation. Since $k_1 + \dots + k_p = p - 1$, there exists an index $j_1 \in \{1,\dots,p\}$ such that $k_{j_1} = 0$. Let
$$(k_1', \dots,k_p') = (k_{j_1+1},\dots,k_p,k_1,\dots,k_{j_1})$$
be the cyclic rotation of $(k_1,\dots,k_p)$, and correspondingly
$$(v_1', \dots,v_p') = (v_{j_1+1},\dots,v_p,v_1,\dots,v_{j_1})$$
be the cyclic rotation of $(v_1,\dots,v_p)$.

Using facts that $\tr(AB) = \tr(BA)$ whenever the dimension of matrices $A,B$ are appropriate, that $Z^{(k_p')} = Z^{(k_{j_1})} = Z^{(0)} = -\mathit{\Pi}$, and that $\mathit{\Pi}P = \mathit{\Pi}$, write
\begin{align*}
-\tr\left(PD^{v_1}Z^{(k_1)}\dots PD^{v_p}Z^{(k_p)}\right) &= -\tr\left(PD^{v_1'}Z^{(k_1')}PD^{v_2'}Z^{(k_2')}\dots PD^{v_{p-1}'}Z^{(k_{p-1}')}PD^{v_p'}Z^{(k_p')}\right)\\
&= \tr\left(PD^{v_1'}Z^{(k_1')}PD^{v_2'}Z^{(k_2')}\dots PD^{v_{p-1}'}Z^{(k_{p-1}')}PD^{v_p'}\mathit{\Pi}\right)\\
&= \tr\left(D^{v_1'}Z^{(k_1')}PD^{v_2'}Z^{(k_2')}\dots PD^{v_{p-1}'}Z^{(k_{p-1}')}PD^{v_p'}\mathit{\Pi}P\right)\\
&= \tr\left(D^{v_1'}Z^{(k_1')}PD^{v_2'}Z^{(k_2')}\dots PD^{v_{p-1}'}Z^{(k_{p-1}')}PD^{v_p'}\mathit{\Pi}\right)\\
&= \langle f, D^{v_1'-1}Z^{(k_1')}PD^{v_2'}Z^{(k_2')}\dots PD^{v_{p-1}'}Z^{(k_{p-1}')}PD^{v_p'-1}f\rangle_\pi.
\end{align*}
Since $k_1'+\dots+k_{p-1}' = p-1 \ge 1$ given $p \ge 2$, there exists an index $j_2 \in \{1,\dots,p-1\}$ such that $k_{j_2}' \ge 1$. From the fact that $Z\mathit{\Pi} = 0$ (\cref{lemma:resolvent}, part (a)), it follows that $Z^{(k_{j_2}')}P = Z^{(k_{j_2}')}(P-\mathit{\Pi})$. Thus,
\begin{align*}
-\tr\left(PD^{v_1}Z^{(k_1)}\dots PD^{v_p}Z^{(k_p)}\right) &= \langle f, D^{v_1'-1}Z^{(k_1')}PD^{v_2'}Z^{(k_2')}\dots Z^{(k_{j_2}')}P \dots PD^{v_{p-1}'}Z^{(k_{p-1}')}PD^{v_p'-1}f\rangle_\pi\\
&= \langle f, D^{v_1'-1}Z^{(k_1')}PD^{v_2'}Z^{(k_2')}\dots Z^{(k_{j_2}')}(P-\mathit{\Pi}) \dots PD^{v_{p-1}'}Z^{(k_{p-1}')}PD^{v_p'-1}f\rangle_\pi\\
&\le \Vert f \Vert_\pi^2 \opnorm{D}_\pi^{(\sum_{j=1}^p v_j') -2} \opnorm{Z}_\pi^{\sum_{j=1}^{p-1} k_j'} \opnorm{P-\mathit{\Pi}}_\pi \opnorm{P}_\pi^{p-2}
\end{align*}
Putting it together with facts that $\Vert f \Vert_\pi^2 = \sigma^2$, $\opnorm{D}_\pi \le c$, $\sum_{j=1}^p v_j' = m$, $\opnorm{Z}_\pi = (1-\rs)^{-1}$ (\cref{lemma:resolvent}, part (d)),  $\sum_{j=1}^{p-1} k_j' = p-1$, $\opnorm{P-\mathit{\Pi}}_\pi = \as$ and $\opnorm{P}_\pi = 1$ yields
$$-\tr\left(PD^{v_1}Z^{(k_1)}\dots PD^{v_p}Z^{(k_p)}\right) \le \sigma^2 c^{m-2} (1-\rs)^{-(p-1)} \as.$$

Combining this result with the following facts
\begin{itemize}
    \item $|\{(v_1,\dots,v_p): v_1+\dots+v_p = m, v_j \ge 1\}| = \binom{m-1}{p-1}$\\
    \item $|\{(k_1,\dots,k_p): k_1+\dots+k_p = p-1, k_j \ge 0\}| = \binom{2p-2}{p-1}$\\
    \item $\sum_{p=1}^m \frac{1}{p}\binom{m-1}{p-1}\binom{2p-2}{p-1} \le 5^{m-2}$ for $m\ge 3$, see \cite[page 856]{lezaud1998chernoff}.\\
    \item $v_1!v_2!\dots v_p! \ge 2^{v_1-1}2^{v_2-1}\dots 2^{v_p-1} = 2^{(v_1+\dots+v_p)-p} = 2^{m-p}$\\
    \item $-1 \le \rs \le 1 \Rightarrow 0 \le \frac{1 - \rs}{2} \le 1$.
\end{itemize}
yield
\begin{equation*}
\begin{split}
\sum_{p=2}^m \beta^{(m)}_p &= \sum_{p=2}^m \frac{1}{p} \sum_{\begin{split}
v_1+\dots+v_p = m, ~v_j \ge 1\\
k_1+\dots+k_p = p-1,~k_j \ge 0\end{split}} \frac{\sigma^2 c^{m-2} (1-\rs)^{-(p-1)}\as}{v_1!v_2!\dots v_p!}\\
&\le 5^{m-2} \times \frac{\sigma^2 c^{m-2} (1-\rs)^{-(p-1)}\as}{2^{m-p}}\\
&= \left(\frac{1-\rs}{2}\right)^{m-p} \frac{\sigma^2\as}{5c}\left(\frac{5c}{1-\rs}\right)^{m-1}\\
&\le \frac{\sigma^2\as}{5c}\left(\frac{5c}{1-\rs}\right)^{m-1}.
\end{split}
\end{equation*}
Further,
$$\beta^{(m)} = \beta^{(m)}_1 + \sum_{p=2}^m \beta^{(m)}_p \le \frac{\pi(f^m)}{m!} + \frac{\sigma^2\as}{5c}\left(\frac{5c}{1-\rs}\right)^{m-1}.$$
The above inequality also holds for $m = 2$, as
$$\beta^{(2)} = \frac{\sigma^2_\text{as}(f)}{2} \le \sigma^2\left(\frac{1}{2} + \frac{\rs}{1-\rs}\right) = \frac{\pi(f^2)}{2!} + \frac{\sigma^2\rs}{5c}\left(\frac{5c}{1-\rs}\right)^{2-1} \le \frac{\pi(f^2)}{2!} + \frac{\sigma^2\as}{5c}\left(\frac{5c}{1-\rs}\right)^{2-1}.$$
Thus,
\begin{align*}
\beta(t) &= \beta^{(0)} + \beta^{(1)}t + \sum_{m = 2}^\infty \beta^{(m)}t^m\\
&\le 1 + 0 + \sum_{m = 2}^\infty \frac{\pi(f^m)t^m}{m!} +  \sum_{m = 2}^\infty \frac{\sigma^2 \as t}{5c} \left(\frac{5ct}{1-\rs}\right)^{m-1}\\
&\le \exp\left(\sum_{m = 2}^\infty \frac{\pi(f^m)t^m}{m!} + \sum_{m = 2}^\infty \frac{\sigma^2 \as t}{5c} \left(\frac{5ct}{1-\rs}\right)^{m-1}\right).
\end{align*}
In the exponent, the first term
$$\sum_{m = 2}^\infty \frac{\pi(f^m)t^m}{m!} \le \sum_{m = 2}^\infty \frac{\pi(f^2)c^{m-2}t^m}{m!} = \frac{\sigma^2}{c^2} \sum_{m = 2}^\infty \frac{c^m t^m}{m!} = \frac{\sigma^2}{c^2}\left(e^{tc} - 1 - tc\right),$$
and the second term, for any $0 \le t < \frac{1-\rs}{5} \times c^{-1} < t_0 = \frac{1 - \rs}{3 - \rs} \times c ^{-1}$,
$$\sum_{m = 2}^\infty \frac{\sigma^2 \as t}{5c} \left(\frac{5ct}{1-\rs}\right)^{m-1} = \frac{\sigma^2 \as t^2}{1- \rs - 5ct},$$
completing the proof.
\end{proof}

\begin{lemma} \label{lemma:infimal convolution}
For $\lambda \in [0,1)$, define for any $0 \le t < (1-\lambda)/5c$
\begin{equation} \label{equation:g}
g_1(t) = \frac{\sigma^2}{c^2}(e^{tc} - 1 - tc), ~~~g_2(t)=\frac{\sigma^2 \lambda t^2}{1-\lambda - 5ct}.
\end{equation}
If $\lambda \in (0,1)$ then
\begin{equation} \label{equation:infimal convolution}
\begin{split}
(g_1+g_2)^*(\epsilon)
&\coloneqq \sup\{t\epsilon - g_1(t) - g_2(t): ~0 \le t < (1-\lambda)/5c\}\\
&\ge \frac{\epsilon^2}{2} \left(\frac{1+\lambda}{1-\lambda} \cdot \sigma^2 + \frac{5c\epsilon}{1-\lambda}\right)^{-1}
\end{split}
\end{equation}
If $\lambda = 0$ then
$$(g_1+g_2)^*(\epsilon) = g_1^*(\epsilon) \ge \frac{\epsilon^2}{2} \left(\sigma^2 +\frac{c\epsilon}{3}\right)^{-1}.$$
\end{lemma}

\begin{proof}[Proof of \cref{lemma:infimal convolution}]
Extend functions $g_1(t)$ and $g_2(t)$ to the domain $(-\infty,+\infty)$ as
\begin{align*}
g_1(t) &= \begin{cases} 0 & \mbox{if } t < 0\\
\frac{\sigma^2}{c^2}\left(e^{tc} - 1 - tc\right) & \mbox{if } t \ge 0
\end{cases}\\
g_2(t) &= \begin{cases} 0 & \mbox{if } t < 0\\
\frac{\sigma^2 \lambda t^2}{1 - \lambda - 5ct} & \mbox{if } 0 \le t < \frac{1-\lambda}{5c}\\
+\infty & \mbox{if } t \ge \frac{1-\lambda}{5c}
\end{cases}
\end{align*}
Both $g_1(t)$ and $g_2(t)$ are closed proper convex functions. Their convex conjugates are given by
\begin{equation} \label{equation:g1star}
g_1^*(\epsilon_1) = \begin{cases}
\frac{\sigma^2}{c^2}h_1\left(\frac{c\epsilon_1}{\sigma^2}\right) & \mbox{if } \epsilon_1 \ge 0,\\
+\infty & \mbox{if } \epsilon_1 <  0,
\end{cases}
\end{equation}
with $h_1(u) = (1+u)\log(1+u) - u$ for $u \ge 0$, and
\begin{equation} \label{equation:g2star}
g_2^*(\epsilon_2) = \begin{cases}
\frac{(1-\lambda)\epsilon_2^2}{2\lambda \sigma^2} h_2\left(\frac{5c\epsilon_2}{\lambda\sigma^2}\right) & \mbox{if } \epsilon_2 \ge 0,\\
+\infty & \mbox{if } \epsilon_2 < 0,
\end{cases}
\end{equation}
with $h_2(u) = (\sqrt{1 + u} + u/2 + 1)^{-1}$.

The convex conjugate $(g_1+g_2)^*$ is still well-defined under this extension. Indeed, $g_1(t) = O(t^2)$ and $g_2(t) = O(t^2)$ as $t \to 0^+$, $t \epsilon - g_1(t) - g_2(t) > 0$ for small enough $t > 0$; and $t \epsilon - g_1(t) - g_2(t) \le 0$ for $t \le 0$. Therefore,
$$(g_1+g_2)^*(\epsilon) \coloneqq \sup_{0 \le t < (1- \lambda)/5c} \{t \epsilon - g_1(t) - g_2(t)\} = \sup_{t \in \mathbb{R}} \{t \epsilon - g_1(t) - g_2(t)\}.$$

By the Moreau-Rockafellar formula \citep[Theorem 16.4]{rockafellar2015convex}, the convex conjugate of $g_1+g_2$ is the infimal convolution of their conjugates $g_1^*$ and $g_2^*$. That is,
$$(g_1+g_2)^*(\epsilon) = \inf \{g_1^*(\epsilon_1) + g_2^*(\epsilon_2):~ \epsilon_1 + \epsilon_2 = \epsilon, ~\epsilon_1 \in \mathbb{R},~\epsilon_2 \in \mathbb{R}\}.$$
Putting \eqref{equation:g1star} and \eqref{equation:g2star} into the last display yields that
$$(g_1+g_2)^*(\epsilon)= \inf \left\{\frac{\sigma^2}{c^2}h_1\left(\frac{c\epsilon_1}{\sigma^2}\right) + \frac{(1-\lambda)\epsilon_2^2}{2\lambda \sigma^2} h_2\left(\frac{5c\epsilon_2}{\lambda\sigma^2}\right): ~\epsilon_1 + \epsilon_2 = \epsilon, ~\epsilon_1 \ge 0, ~\epsilon_2 \ge 0\right\}.$$
Bounding $h_1(u) \ge \frac{u^2}{2(1+u/3)}$ and $h_2(u) \ge \frac{1}{2 + u}$ for $u \ge 0$ yields that
$$(g_1+g_2)^*(\epsilon) \ge \inf \left\{ \frac{\epsilon_1^2}{2\left(\sigma^2 + \frac{c\epsilon_1}{3}\right)}+ \frac{\epsilon_2^2}{2\left( \frac{2\lambda}{1-\lambda} \sigma^2 + \frac{5c\epsilon_2}{1-\lambda}\right)}: ~\epsilon_1 + \epsilon_2 = \epsilon, ~\epsilon_1 \ge 0, ~\epsilon_2 \ge 0\right\}$$
Using the fact that $\epsilon_1^2/a + \epsilon_2^2/b \ge (\epsilon_1+\epsilon_2)^2/(a+b)$ for any non-negative $\epsilon_1,\epsilon_2$ and positive $a,b$ yields
\begin{align*}
(g_1+g_2)^*(\epsilon) &\ge \inf \left\{ \frac{(\epsilon_1+\epsilon_2)^2}{2\left(\sigma^2 + \frac{c\epsilon_1}{3}\right) + 2\left( \frac{2\lambda}{1-\lambda} \sigma^2 + \frac{5c\epsilon_2}{1-\lambda}\right)}: ~\epsilon_1 + \epsilon_2 = \epsilon, ~\epsilon_1 \ge 0, ~\epsilon_2 \ge 0\right\}\\
&= \frac{\epsilon^2}{2\left( \frac{1+\lambda}{1-\lambda} \sigma^2 + \frac{5c\epsilon}{1-\lambda}\right)}.
\end{align*}
For the case of $\lambda = 0$, we merely need to lower bound $h_1(u) \ge \frac{u^2}{2(1+u/3)}$ as $g_1^*(\epsilon) = (g_1+g_2)^*(\epsilon)$.
\end{proof}

\begin{proof}[Proof of \cref{thm1}]
With \cref{lemma:norm-product1} in hand, it suffices to show that, for any function $f$ with $|f| \le c$, $\pi(f) = 0$ and $\pi(f^2) = \sigma^2$,
\begin{equation} \label{equation:operator norm}
\log \opnorm{E^{tf/2}\widehat{P}_\lambda E^{tf/2}}_\pi \le g(t; \as, c, \sigma^2) := \frac{\sigma^2}{c^2} (e^{tc} - 1 - tc) + \frac{\sigma^2 \as t^2}{1- \as - 5ct}
\end{equation}
To this end, we discretize $f$ as $f_k$ in the following way. Cut the value range $[-c, +c]$ of function $f$ into $2k$ bins of length $c/k$. Let $m_k(x)$ be the midpoint of $i$-th bin if $f(x)$ belongs to the $i$-th bin, and let $f_k = m_k - \pi(m_k)$. It is not hard to verify that $f_k$ uniformly converges to $f$. From \cref{lemma:discretization}, it follows that
$$\opnorm{E^{tf/2}\widehat{P}_\as E^{tf/2}}_\pi = \lim_{k \to \infty} \opnorm{E^{tf_k/2}\widehat{P}_\as E^{tf_k/2}}_\pi.$$
For each $E^{tf_k/2}\widehat{P}_\as E^{tf_k/2}$, let $\{\widehat{X}_i\}_{i \ge 1}$ be a Markov chain of $\widehat{P}_\as$ and $\widehat{Y}_i = f_k(\widehat{X}_i)$. Then $\{Y_i\}_{i \ge 1}$ is a Markov chain in a finite state space $\mc{Y}_k = f_k(\mc{X})$, as $f_k$ takes at most $2k$ different values. Let $\widehat{Q}_\as$ and $\mu_k$ denote the transition kernel (transition probability matrix) and invariant distribution (probability vector) of the Markov chain $\{\widehat{Y}_i\}_{i \ge 1}$, and let $E^{t\mathcal{Y}_k/2}$ be the diagonal matrix of $e^{ty/2}$ for $y \in \mc{Y}_k$. By \cref{lemma:counterpart},
$$\opnorm{E^{tf_k/2}\widehat{P}_\as E^{tf_k/2}}_\pi = \opnorm{E^{t\mathcal{Y}_k/2}\widehat{Q}_\as E^{t\mathcal{Y}_k/2}}_{\mu_k}.$$
Note that
\begin{align*}
\sum_{y \in \mc{Y}_k} \mu_k y &= \sum_{y \in \mc{Y}_k} \pi(\{x: f_k(x) = y\}) y = \pi(f_k) = 0\\ 
\sum_{y \in \mc{Y}_k} \mu_k y^2 &= \sum_{y \in \mc{Y}_k} \pi(\{x: f_k(x) = y\}) y^2 = \pi(f_k^2)
\end{align*}
and that $\widehat{Q}_\as$ is a reversible, irreducible transition probability matrix on finite state space $\mathcal{Y}_k$ with both absolute and right spectral gaps being $1-\as$. It follows from \cref{lemma:finite-state-space} that
$$\log \opnorm{E^{t\mathcal{Y}_k/2}\widehat{Q}_\as E^{t\mathcal{Y}_k/2}}_{\mu_k} \le g\left(t; \as, c, \pi(f_k^2)\right).$$
Collecting these pieces together and letting $\pi(f_k^2) \to \sigma^2$ as $k \to \infty$ yield \eqref{equation:operator norm}.
\end{proof}

\begin{proof}[Proof of \cref{thm2}]
Just substitute \cref{lemma:norm-product1} with \cref{lemma:norm-product2} in the proof of \cref{thm1}.
\end{proof}

\begin{proof}[Proof of \cref{thm3}]
Since $\mathbb{E}_{\nu}\left[\cdot\right] = \mathbb{E}_{\pi}\left[\cdot\left(\frac{d\nu}{d\pi}\right)\right]$, an elementary use of the H\"{o}lder's inequality 
yields
$$\mathbb{E}_{\pi}\left[e^{t\sum_{i=1}^{n} f(X_i)}\right] \le e^{n g(t)} \Rightarrow \mathbb{E}_{\nu}\left[e^{t\sum_{i=1}^{n} f(X_i)}\right] \le \left\Vert \frac{d\nu}{d\pi} \right\Vert_{\pi,p} e^{n g(qt)/q}.$$
The tail probability bound follows from the fact that the convex conjugate of $g_q(t) = g(qt)/q$ is given by $g^*_q(\epsilon) = g^*(\epsilon)/q$ \citep[Page 95]{boyd2004convex}.
\end{proof}

\begin{proof}[Proof of \cref{thm4}]
The proof is similar to that of \citet*[Theorem 4]{fan2021hoeffding}, and thus is omitted.
\end{proof}

\section{Applications}\label{sec5}

\scolor{
\subsection{MCMC integral estimation with non-asymptotic confidence intervals}

MCMC methods are widely used for numerical approximation of multi-dimensional integrals in Bayesian statistics, computational physics and computational biology \citep{gilks1995markov}. To showcase the advantages of our new Bernstein-type inequalities, we derive an  MCMC estimator with a non-asymptotic confidence interval for a multi-dimensional integral and compare this confidence interval with those derived by previous Bernstein-type and Hoeffding-type inequalities. 

Consider a typical example of the MCMC estimation for an integral in $\mathbb{R}^3$
$$\pi(f) = \int_0^1 \int_0^1 \int_0^1 f(x,y,z)\pi(x,y,z)dxdydz,$$
where
$$f(x,y,z) = \sqrt{\frac{x^2+y^2+z^2}{3}}$$
and
$$\pi(x,y,z) \propto \exp\left(\sin\left(\frac{\pi xyz}{2}\right)\right).$$

A reversible Markov chain $\{(X_i, Y_i, Z_i)\}_{i \ge 1}$ is generated by an independent Metropolis-Hastings (IMH) algorithm with the uniform distribution
$$\nu(x,y,z) = 1,~~~\forall~(x,y,z) \in [0,1]^3$$
as the initial distribution and the proposal kernel. The details are collected in \cref{alg1}. Note that the normalizing constant of the probability density function $\pi(x,y,z)$ is not needed for the implementation of \cref{alg1}, as only the ratio of $\pi(x',y',z')/\pi(x,y,z)$ is used in the calculation of the acceptance probability.

\begin{algorithm}[t] 
\caption{Independent Metropolis-Hastings (IMH) Algorithm} 
\label{alg1} 
\begin{algorithmic} 
\State Pick an initial state $(x_0,y_0,z_0) \sim \nu(x,y,z)$.
\For{$i = 0, 1, 2, \dots$}
    \State Generate a random candidate state $(x', y', z') \sim \nu(x,y,z)$.
    \State Calculate the acceptance probability $A(x_i, y_i, z_i, x', y', z') = \min\left\{1, \frac{\pi(x',y',z')\nu(x_i,y_i,z_i)}{\pi(x_i,y_i,z_i)\nu(x',y',z')}\right\}$
    \State Generate a uniform random number $u \in [0,1]$.
    \If{$u \le A(x_i,y_i,z_i, x', y', z')$}
        \State accept the new state $(x',y',z')$ and set $(x_{i+1}, y_{i+1}, z_{i+1}) = (x',y',z')$.
    \Else
        \State reject the new state $(x',y',z')$ and set $(x_{i+1}, y_{i+1}, z_{i+1}) = (x_i, y_i, z_i)$.
    \EndIf
\EndFor
\end{algorithmic}
\end{algorithm}

An MCMC estimator for the desired integral $\pi(f)$ is given by
$$\hat{f}_n = \frac{1}{n}\sum_{i=0}^{n-1} f(X_{ik+1},Y_{ik+1}, Z_{ik+1})$$
for some positive integer $k$. Using evenly-spaced subsamples of a Markov chain is a popular trick in the practice of MCMC \citep{geyer1992practical}. These sub-samples $\{(X_{ik+1},Y_{ik+1},Z_{ik+1})\}_{i \ge 0}$ form a Markov chain with a transition kernel $P^k$, which is the $k$-th iterate of the transition kernel $P$ of the original Markov chain. It is clear that $\rs(P^k) = \rs(P)^k$ by the spectral mapping theorem.

Applying the Bernstein-type inequality in \cref{thm3} yields
\begin{align*}
\P_\nu\left(| \hat{f}_n - \pi(f)| \ge \epsilon\right) \le 2 \sup_{(x,y,z) \in [0,1]^3} \frac{\nu(x,y,z)}{\pi(x,y,z)} \times \exp\left(-\frac{n\epsilon^2/2}{\frac{1+\lambda^k}{1-\lambda^k}\sigma^2 + \frac{5\epsilon}{1-\lambda^k}}\right).
\end{align*}
Let $C$ be the normalizing constant of $\pi$, namely,
$$C \coloneqq \int_0^1\int_0^1\int_0^1 \exp\left(\sin\left(\frac{\pi xyz}{2}\right)\right)dxdydz$$
then
$$\sup_{(x,y,z) \in [0,1]^3} \frac{\nu(x,y,z)}{\pi(x,y,z)} = \frac{\nu(0,0,0)}{\pi(0,0,0)} = C.$$

It follows that, with probability at least $1-2\delta$, $\pi(f)$ is within the interval
$$\hat{f}_n \pm \sqrt{\frac{1+\lambda^k}{1-\lambda^k}} \cdot \sigma \cdot \sqrt{\frac{2\log(C/\delta)}{n}} + \frac{5}{1-\lambda^k} \cdot \frac{2\log(C/\delta)}{n}.$$

For the unknown constant $C$, we substitute it with its upper bound $e$. For the unknown constant $\sigma^2 \coloneqq \text{Var}_\pi(f)$, we substitute it with its sample estimate. For the unknown constant $\lambda$, we substitute it with its upper bound $1 - e^{-1}$. This upper bound is deduced by using \citet[Theorem 2.2]{atchade2007geometric}. Write
$$
1 - \as(P) = 1 - \rs(P) = 1- \lambda = \inf_{(x,y,z) \in [0,1]^3} \bar{A}(x,y,z),
$$
where $\bar{A}(x,y,z)$ is the overall acceptance probability at state $(x,y,z)$
\begin{align*}
\bar{A}(x,y,z) &= \int_0^1\int_0^1\int_0^1 A(x,y,z,x',y',z') \nu(x',y',z')dx'dy'dz'\\
&= \int_0^1\int_0^1\int_0^1 \min\left\{1, \frac{\pi(x',y',z')\nu(x,y,z)}{\pi(x,y,z)\nu(x',y',z')}\right\}  \nu(x',y',z')dx'dy'dz'
\end{align*}
Apparently, $(x_\star, y_\star,z_\star) = (1,1,1)$ minimizes the density ratio $\frac{\nu(x,y,z)}{\pi(x,y,z)}$. Thus
\begin{align*}
\inf_{(x,y,z) \in [0,1]^3} \bar{A}(x,y,z) &= \int_0^1 \int_0^1 \int_0^1 \min \left\{1, \frac{\pi(x',y',z') / \nu(x',y',z')}{ \pi(x_\star,y_\star,z_\star) / \nu(x_\star, y_\star,z_\star)}\right\} \nu(x',y',z')dx'dy'dz'\\
&= \int_0^1 \int_0^1 \int_0^1 \frac{\pi(x',y',z') / \nu(x',y',z')}{ \pi(x_\star,y_\star,z_\star) / \nu(x_\star, y_\star, z_\star)} \nu(x',y',z')dx'dy'dz'\\
&= \frac{\nu(x_\star, y_\star, z_\star)}{\pi(x_\star, y_\star, z_\star)} = \inf_{(x,y,z) \in [0,1]^3} \frac{\nu(x,y,z)}{\pi(x,y,z)} = Ce^{-1} \ge e^{-1}.
\end{align*}

With $k=10$, $\delta = 0.005$ and $n = 100000$, experimental results show that $\pi(f)$ is within the interval $0.5829 \pm 0.0052$ with at least 99\% probability. \cref{tab} compares this confidence interval with those derived by other Bernstein-type and Hoeffding-type inequalities for Markov chains. The confidence interval derived by our Bernstein-type inequalities is at least $30\%$ tighter than those derived by previous Bernstein-type and Hoeffding-type inequalities.

\begin{table}[h!]
    \centering
    \begin{tabular}{l c c}
         \hline
         &  Inequality Type & Half Length of CI\\
         \hline
         \cref{thm3} & Bernstein & 0.0052\\
         \citet[Equation 3.20]{paulin2015concentration} & Bernstein &  0.0068\\
         \citet[Equation 3.21]{paulin2015concentration} & Bernstein & 0.0071\\
         \citet[Equation 3.22]{paulin2015concentration} & Bernstein & 0.0104\\ 
         \citet*[Theorem2]{fan2021hoeffding} & Hoeffding & 0.0080\\
         \hline
    \end{tabular}
    \caption{Comparing confidence intervals (CIs) derived by selective concentration inequalities for Markov chains}
    \label{tab}
\end{table}

As a final remark, we note that, although a change-of-variable trick $t = xyz$ and some numerical integration method can approximate
$$C = \int_0^1 \exp\left(\sin\left(\frac{\pi t}{2}\right)\right) \frac{\log^2(t)}{2}dt \approx 1.2335$$ in this typical example and further $\lambda = 1-Ce^{-1} \approx 0.5462$, the normalizing constant $C$ is usually difficult to compute in general cases. We use the upper bounds $\lambda \le 1-e^{-1}\approx 0.6321$ and $C\leq e\approx 2.7183$ rather than their approximate values $\lambda \approx 0.5462$ and $C\approx 1.2335$ to construct the confidence intervals in \cref{tab}. Nonetheless, the approach for the IMH algorithm showcased in this section is applicable whenever the density ratio $\nu/\pi$ has known upper and lower bounds.}

\subsection{Robust mean estimation under Markov dependence}\label{sec:robust}

\scolor{In our second example, w}e consider to robustly estimate the mean of a sequence of Markov-dependent samples. 
Suppose the data $y_1,\ldots, y_n\in \RR$ are generated according to 
\#\label{model:mean}
y_i = \mu^* + \varepsilon_i, 
\#
where $\mu^*\in \RR$ is the underlying mean and $\{\varepsilon_i: 1\leq i\leq n\}$ with $\mathbb{E}\, \varepsilon_i =0$ is a stationary and general Markov chain  with invariant measure $\pi$ and absolute spectral gap $1-\as>0$. \scolor{With a slight overload of notation, we assume $\pi$ has only bounded second moment such that 
\$
\sigma^2\coloneqq \int x^2 \pi(d x) <\infty.
\$

To estimate the mean $\mu^*$, the sample mean estimator $\sum_{i=1}^n y_i/n$ is known to achieve at best a polynomial-type non-asymptotic confidence width even in the independent case \citep{catoni2012challenging}. That is, for i.i.d. $y_i$, there is some distribution $F$ for $\varepsilon_i$ with mean $0$ and variance $\sigma^2$ such that
\$
\PP\left(\left|\sum_{i=1}^n \frac{y_i}{n} -\mu^*\right|\leq c \sigma\sqrt{\frac{1}{n}\cdot\frac{1}{\delta}}\right) > 1-2\delta,
\$ 
for some constant $c$.  Intuitively, this means that the sample mean does not concentrate on the true mean fast enough when errors have only finite variances. \cite{catoni2012challenging} made an important step towards estimating the mean with faster concentration. Specifically, Catoni constructs a robust mean estimator $\widehat{\mu}_\tau$, depending on some tuning parameter $\tau$, that deviates from the true mean $\mu$ logarithmically in $1/\delta$, that is
\# \label{catoni.bound}
	 \PP\left(   |  \widehat{\mu}_\tau - \mu^* | \leq  c\sigma \sqrt{\frac{1}{n}\cdot \log\left(\frac{1}{\delta}\right)}\right) > 1- 2\delta,
\#
for some constant $c$. 
}

\scolor{
However, no such results exist for Markov-dependent data with only second moment.  Our goal is to estimate $\mu^*$ in \eqref{model:mean} with a similar non-asymptotic guarantee as in \eqref{catoni.bound}.} Similar to \cite{sun2020adaptive}, we will make use of the Huber loss: $\ell_\tau(x)= x^2/2$ if $|x|\leq \tau$; and $\ell_\tau(x)= \tau|x|-\tau^2/2$, elsewhere. We consider the following estimator
\$
\widehat\mu _\tau =\arg\min_\mu \left\{\frac{1}{n}  \sum_{i=1}^n \ell_\tau(y_i-\mu) =: L_n(\mu)\right\}. 
\$

\scolor{
 We need the following curvature lemma, whose proof is provided later in this subsection. 
\begin{lemma}\label{lemma:sc}
Suppose that $\tau\geq c(r + \sigma)$ and $n\geq c(1-\lambda)^{-1}(1+\lambda) \log (1/\delta)$ for some absolute constant $c$. 
With probability at least $1-\delta$, we have 
\$
\kappa_-(\mu^*, r)= \inf_{\mu} \left\{ \frac{\left\langle  \nabla L_n(\mu)-\nabla L_n(\mu^*) , \mu -\mu^* \right\rangle}{|\mu-\mu^*|^2}:  \mu \in \BB_r(\mu^*)\right\}\geq 1/2.
\$
\end{lemma}
}

\scolor{
The above lemma indicates that the loss function $L_n(\mu)$ is locally one-point strongly convex with respect to (w.r.t.) a single point $\mu^*$.   Intuitively, this lemma should hold with high probability since the Huber loss is quadratic in the center and linear at tails. In the independent case,  \cite{sun2020adaptive} showed a similar result  holds with high probability in the independent case. With this lemma,  we are ready to derive the non-asymptotic properties of $\widehat\mu_\tau$. 
}

\begin{theorem}\label{thm:robust}
Take $\tau = c_0 \sigma \sqrt{n/\log(1/\delta)}$. Suppose  $n\geq c\big( \sqrt{\alpha_1(\lambda)}\vee \alpha_2(\lambda) \big) \log(1/\delta)$ for some constant $c$, where $\alpha_1(\lambda)$ and $\alpha_2(\lambda)$ are defined in \cref{thm1}, that is 
$$\alpha_1(\as) = \frac{1+\as}{1-\as},~~~\alpha_2(\as) = \begin{cases}
1/3 &\mbox{if}~\as = 0,\\
5/(1-\as) &\mbox{if}~\as \in (0,1).
\end{cases}$$
 Then with probability at least $1-3\delta$ 
\$
\left|\widehat \mu_\tau - \mu^*\right|\leq C \sigma\left(\sqrt{\alpha_1(\lambda)}\vee\alpha_2(\lambda)\right) \sqrt{\frac{\log(1/\delta)}{n}},
\$
for some constant $C$. 
\end{theorem}

The theorem above is a direct consequence of our newly-derived Bernstein-type inequalities, while using the Hoeffding-type inequality for general Markov chains by \citet*{fan2021hoeffding} would give a significantly weaker, in fact non-diminishing, upper bound. Applying the Hoeffdinig-type inequality in the proof of \cref{thm:robust}, we would get instead, with probability at least $1-3\delta$,
\$
\left|\widehat \mu_\tau - \mu^*\right|\leq C \sigma \sqrt{\alpha_1(\lambda)}.
\$

\section{Discussions} \label{sec6}

This paper extends Lezaud's method to derive Bernstein-type inequalities for general-state-space Markov chains. The main challenge addressed is the analysis of the norm of an infinite-dimensional operator. We overcome this difficulty by employing a Leon-Perron operator as the extremal counterpart to the original Markov operator, and by using function discretization techniques to reduce the problem to a finite state space. Both the L\'eon-Perron operator and the discretization technique were originally proposed for Hoeffding-type inequalities of Markov chains \citep*{leon2004optimal,miasojedow2014hoeffding,fan2021hoeffding}. Interestingly, they can help the derivation of Bernstein-type inequalities.

In finite state spaces, we refine Lezaud's analysis of the second-largest eigenvalue of the perturbed transition probability matrix, resulting in tighter moment-generating function (mgf) bounds for Bernstein-type inequalities. In these mgf bounds, the effects of independence and Markov dependence are clearly separated as two distinct convex functions $g_1(t)$ and $g_2(t)$. The variance proxies for these mgf bounds are optimal, achieving the lower bound of variance proxies specified by the Central Limit Theorem for Markov chains.

A further challenge emerges when converting the mgf bound into a convenient tail probability bound, primarily due to the absence of an explicit expression for the convex conjugate of $g_1(t) + g_2(t)$. To overcome this, we utilize techniques from convex analysis theory to establish a lower bound. This approach maintains the sharpness of the variance proxy.

We demonstrate the applicability of our newly derived Bernstein-type inequalities with two examples in statistics and machine learning. We anticipate that these inequalities will have additional applications in the future.

\appendix

\vspace{30pt}
\noindent{\bf \LARGE Appendix}

\section{Proofs for Section \ref{sec5}}

This section presents the proofs to Lemma \ref{lemma:sc} and Theorem \ref{thm:robust}. We start with the proof of Lemma \ref{lemma:sc}. 

\begin{proof}[Proof of \cref{lemma:sc}]
For notational simplicity, let $f_i(\mu) = \ell_\tau(y_i -\mu)$ and 
\$
\psi_\tau(x) \coloneqq \nabla \ell_\tau(x) = \sign(x)(|x|\wedge \tau ). 
\$
By convexity,  $\langle \nabla f_i(\mu) - \nabla f_i(\mu^*), \mu - \mu^* \rangle \geq 0$ for any $i$,  and hence
\$
D(\mu) & \coloneqq \langle \nabla L_n(\mu) - \nabla L_n(\mu^*), \mu-\mu^* \rangle \\
& = \frac{1}{n}\sum_{i=1}^n \left\langle \psi_\tau (\varepsilon_i) - \psi_\tau (y_i - \mu), \, \mu - \mu^* \right\rangle \\
& \geq \frac{1}{n}\sum_{i=1}^n \left\langle \psi_\tau (\varepsilon_i) - \psi_\tau (y_i - \mu), \, \mu - \mu^* \right\rangle  1_{\cE_i},
\$
where  $1_{\cE_i}$ is the indicator function of the event
\$
\cE_i = \big\{ |\varepsilon_i|_2 \leq \tau/2 \big\} \cap \big\{|\mu - \mu^*| \leq  r\leq \tau/2 \big\} = \big\{ |\varepsilon_i|_2 \leq \tau/2 \big\}
\$
where the second equality automatically verifies as long as $\tau\geq 2r$. 
On event $\cE_i$,  observe that $\psi_\tau(\varepsilon_i) = \varepsilon_i$ and $\psi_\tau (y_i - \mu) = y_i - \mu$ for all $\mu \in \BB_r(\mu^*)$.   Consequently,
\$
D(\mu) 
&\geq \frac{1}{n}\sum_{i=1}^n \langle \mu -\mu^*    , \mu -\mu^*   \rangle \mathbbm{1}_{\cE_i}\\
&= |\mu -\mu^*|^2 \cdot \frac{1}{n}\sum_{i=1}^n  \mathbbm{1}_{\cE_i}
\$
for any $\mu \in  \BB_r(\mu^*).$
Thus to prove the lemma, it suffices to lower bound
$
\frac{1}{n}\sum_{i=1}^n  \mathbbm{1}_{\cE_i}. 
$
By the Hoeffding-type inequality for general Markov chains \citep*[Theorem 3]{fan2021hoeffding}, we have with probability at least $1-\delta$
\$
\frac{1}{n}\sum_{i=1}^n  \mathbbm{1}_{\cE_i} - \EE 1_{\cE_i} 
&\geq - \sqrt{\frac{1+\lambda}{1-\lambda} \cdot \frac{\log(1/\delta)}{2n}}. 
\$
 By Markov's inequality,
\$
1\geq \PP\big(  |\varepsilon_i | \leq \tau/2 \big) \geq 1 -  \frac{4\sigma^2}{\tau^2} \geq \frac{3}{4}. 
\$
Now because
\$
n\geq 32 \, \frac{1+\lambda}{1-\lambda}\log(1/\delta),
\$
thus it holds that
\$
\frac{1}{n}\sum_{i=1}^n 1_{\cE_i} 
\geq \EE 1_{\cE_i} - \sqrt{\frac{1+\lambda}{1-\lambda}\cdot \frac{\log(1/\delta)}{2n}}\geq \frac{1}{2}. 
\$
This finishes the proof.
\end{proof}

\begin{proof}[Proof of \cref{thm:robust}]
For simplicity of notation, write $\mathbb{E}_\pi$ as $\mathbb{E}$ and  $\mathbb{P}_\pi$ as $\mathbb{P}$.
\scolor{
By Lemma \ref{lemma:sc}, we have, if $\tau\geq c(r + \sigma)$ and $n\geq c(1-\lambda)^{-1}(1+\lambda) \log (1/\delta)$ for some absolute constant $c$, then with probability at least $1-\delta$
\$
\kappa_-(\mu^*, r)= \inf_{\mu} \left\{ \frac{\left\langle  \nabla L_n(\mu)-\nabla L_n(\mu^*) , \mu -\mu^* \right\rangle}{|\mu-\mu^*|^2}:  \mu \in \BB_r(\mu^*)\right\}\geq 1/2.
\$
}

From now on, we implicitly condition on the above probability event $\mathcal{E}$. Take $r$ such that 
\$
r=4\sigma\left(4c_0 \alpha_2(\lambda)+ \sqrt{2\alpha_1(\lambda)}+c_0\right)  \sqrt{\frac{\log(1/\delta)}{n}}  
\$
and assume that $|\widehat\mu_\tau -\mu^*|\leq r$. Then
\$
\frac{1}{2}|\widehat\mu_\tau-\mu^*|^2\leq \left| \nabla L_n(\widehat\mu)- \nabla L_n(\mu^*)\right| \cdot\left |\widehat\mu_\tau-\mu^*\right|,
\$
or equivalently
\$
\frac{1}{2}|\widehat\mu_\tau-\mu^*|
&\leq \left| \nabla L_n(\widehat\mu_\tau)- \nabla L_n(\mu^*)\right|
\leq \left|  \nabla L_n(\mu^*) -\mathbb{E} \nabla L_n(\mu^*)\right| +   \left|  \mathbb{E} \nabla L_n(\mu^*)\right|,
\$
where we use the first order optimality condition that $\nabla L_n(\widehat\mu) =0$  because $\widehat\mu_\tau$ is a minimum of $L_n(\mu)$. Both terms in the right hand side of the above inequality involve the gradient of the loss function, which  can be written as
\#\label{eq:grad}
\nabla L_n(\mu^*) = \frac{1}{n}\sum_{i=1}^n {\rm sign}(\varepsilon_i)\cdot\left(|\varepsilon_i| \wedge \tau\right),
\#
where ${\rm sign}(\cdot)$ is the sign function. We start with the second term. Because $\mathbb{E}(\varepsilon_i) =0$, we have  
\$
\left|\mathbb{E} \nabla L_n(\mu^*)\right| 
&= \frac{1}{n}\left|\sum_{i=1}^n \left\{\mathbb{E} \left( {\rm sign}(\varepsilon_i)\cdot (|\varepsilon_i|\wedge \tau)\right) - \varepsilon_i\right\}  \right| \\
&= \frac{1}{n}\left|\sum_{i=1}^n \mathbb{E} \left(\varepsilon_i\cdot 1(|\varepsilon_i|>\tau) \right) \right| \\
&\leq \frac{1}{n}\sum_{i=1}^n \frac{\mathbb{E}\left(\varepsilon_i^2\cdot 1 (|\varepsilon_i|> \tau)\right)}{\tau}  \leq   \frac{\sigma^2}{\tau},
\$
where the last second inequality is due to Markov inequality. We then bound the  first term. Rewrite $\nabla L_n(\mu^*) - \mathbb{E} L_n(\mu^*)$ as 
\#\label{eq:grad}
\nabla L_n(\mu^*) - \mathbb{E} L_n(\mu^*) \nonumber
&= \frac{1}{n}\sum_{i=1}^n {\rm sign}(\varepsilon_i)\cdot\left(|\varepsilon_i| \wedge \tau\right) - \frac{1}{n}\sum_{i=1}^n \mathbb{E}\left\{{\rm sign}(\varepsilon_i)\cdot\left(|\varepsilon_i| \wedge \tau\right)\right\}\nonumber\\
&= \frac{1}{n}\sum_{i=1}^n f(\varepsilon_i) 
\#
where $f(\varepsilon_i)=   {\rm sign}(\varepsilon_i) \cdot(|\varepsilon_i|\wedge \tau) - \mathbb{E}\left\{ {\rm sign}(\varepsilon_i) \cdot(|\varepsilon_i|\wedge \tau)\right\}: \RR \rightarrow [-\tau-\frac{\sigma^2}{\tau}, \tau+\frac{\sigma^2}{\tau}]$. We have $\sum_{i=1}^n\pi(f^2)/n \leq \sigma^2$. Thus applying \cref{thm1} with $c= \tau + \frac{\sigma^2}{\tau}$ and $f_i=f$, we obtain for any $t > 0$
\begin{equation} \label{eq:robust_tail}
\mathbb{P}\left(\frac{1}{n}\sum_{i=1}^n f(\varepsilon_i)>t\right) \le \exp{\left(-\frac{nt^2/2}{\alpha_1(\as) \cdot \sigma^2 + \alpha_2(\as) \cdot ct}\right)},
\end{equation}
where $\alpha_1(\lambda)$ and $\alpha_2(\lambda)$ are defined in \cref{thm1}. Taking $\tilde t =\frac{nt^2/2}{\alpha_1(\as) \cdot \sigma^2 + \alpha_2(\as) \cdot ct}$ in \eqref{eq:robust_tail} and using a symmetry argument, we obtain for any $\tilde t >0$
\$
\mathbb{P}\left(\left|\frac{1}{n}\sum_{i=1}^n f(\varepsilon_i)\right|> \frac{2\alpha_2(\lambda) c \tilde t}{n} + \sqrt{\frac{2\alpha_1(\lambda)\sigma^2\tilde t}{n}}\right) \le 2\exp (-\tilde t).
\$
 Taking $\delta= \exp(-\tilde t)$ in the above inequality and $\tau = c_0 \sigma \sqrt{n/\log(1/\delta)}\geq\sigma$, we obtain with probability at least $1-2\delta$
\$
\left|\frac{1}{n}\sum_{i=1}^n f(\varepsilon_i)\right|
&\leq  \frac{2\alpha_2(\lambda) c \log(1/\delta)}{n} + \sqrt{\frac{2\alpha_1(\lambda)\sigma^2\log(1/\delta)}{n}}\\
&\leq \frac{4\tau \alpha_2(\lambda) \log(1/\delta)}{n} + \sqrt{2}\sigma \sqrt{\frac{\alpha_1(\lambda)  \log(1/\delta)}{n}}\\
&\leq  \left(4c_0 \alpha_2(\lambda)+ \sqrt{2\alpha_1(\lambda)}\right) \sigma \sqrt{\frac{\log(1/\delta)}{n}}\lesssim \sigma \left(\alpha_2(\lambda)+\sqrt{\alpha_1(\lambda)}\right) \sqrt{\frac{\log(1/\delta)}{n}}.
\$
Combining upper bounds for both terms, we acquire with probability at least $1-2\delta$
\$
|\widehat\mu_\tau-\mu^*|\leq 2\sigma\left(4c_0\alpha_2(\lambda)+\sqrt{2\alpha_1(\lambda)}+c_0\right)\sqrt{\frac{\log(1/\delta)}{n}}. 
\$
\scolor{Without conditioning on the probability event $\cE$, we obtain with probability at least $1-3\delta$ the above inequality holds. }

We then show that  $|\widehat \mu_\tau-\mu^*|\leq r$ must hold. If not, we shall construct an intermediate solution between $\mu^*$ and $\widehat\mu_\tau$, denoted by $\mu_{\tau,\eta} =\mu^*+ \eta (\widehat\mu_\tau-\mu^*)$, such that $| \mu_{\tau, \eta} -\mu^* |= r$. Specifically, we can  choose some $\eta \in(0,1)$ so that  $|  \mu_{\tau,\eta} -\mu^* |=r$. We then repeat  the above calculation and  obtain 
\$
|  \mu_{\tau,\eta} -\mu^* |\leq  2\sigma\left(4c_0 \alpha_2(\lambda)+ \sqrt{2\alpha_1(\lambda)}+c_0\right)  \sqrt{\frac{\log(1/\delta)}{n}}     <r. 
\$
which is a  contradiction. Thus it must  hold that $|\widehat \mu_\tau-\mu^*|\leq r$. Finally the conditions in the theorem verify the conditions used in the proof and   this completes the proof. 

\end{proof}



\bibliographystyle{apalike}
\bibliography{ref}

\begin{thebibliography}{}

\bibitem[Adamczak and Bednorz, 2015]{adamczak2015exponential}
Adamczak, R. and Bednorz, W. (2015).
\newblock Exponential concentration inequalities for additive functionals of
  {M}arkov chains.
\newblock {\em ESAIM: Probability and Statistics}, 19:440--481.

\bibitem[Aldous, 1983]{aldous1983random}
Aldous, D. (1983).
\newblock Random walks on finite groups with slowly increasing orders.
\newblock {\em Annals of Probability}, 11(4):1063--1073.

\bibitem[Aldous, 1981]{aldous1981representations}
Aldous, D.~J. (1981).
\newblock Representations for partially exchangeable arrays of random
  variables.
\newblock {\em Journal of Multivariate Analysis}, 11(4):581--598.

\bibitem[Atchad{\'e} and Perron, 2007]{atchade2007geometric}
Atchad{\'e}, Y.~F. and Perron, F. (2007).
\newblock On the geometric ergodicity of metropolis-hastings algorithms.
\newblock {\em Statistics}, 41(1):77--84.

\bibitem[Baum et~al., 1970]{baum1970maximization}
Baum, L.~E., Petrie, T., Soules, G., and Weiss, N. (1970).
\newblock A maximization technique occurring in the statistical analysis of
  probabilistic functions of markov chains.
\newblock {\em The Annals of Mathematical Statistics}, 41(1):164--171.

\bibitem[Bennett, 1962]{bennett1962probability}
Bennett, G. (1962).
\newblock Probability inequalities for the sum of independent random variables.
\newblock {\em Journal of the American Statistical Association},
  57(297):33--45.

\bibitem[Bernstein, 1946]{bernstein1946theory}
Bernstein, S.~N. (1946).
\newblock {\em Theory of Probability (in Russian)}.
\newblock Moscow.

\bibitem[Bishop, 2006]{bishop2006pattern}
Bishop, C.~M. (2006).
\newblock {\em Pattern Recognition and Machine Learning}.
\newblock Springer.

\bibitem[Bottou, 2010]{bottou2010large}
Bottou, L. (2010).
\newblock Large-scale machine learning with stochastic gradient descent.
\newblock In {\em Proceedings of the 18th International Conference on
  Computational Statistics}.

\bibitem[Boucheron et~al., 2013]{boucheron2013concentration}
Boucheron, S., Lugosi, G., and Massart, P. (2013).
\newblock {\em Concentration Inequalities: A Nonasymptotic Theory of
  Independence}.
\newblock Oxford University Press.

\bibitem[Boyd and Vandenberghe, 2004]{boyd2004convex}
Boyd, S. and Vandenberghe, L. (2004).
\newblock {\em Convex optimization}.
\newblock Cambridge university press.

\bibitem[Breiman, 1992]{breiman1992probability}
Breiman, L. (1992).
\newblock {\em Probability}.
\newblock Society for Industrial and Applied Mathematics.

\bibitem[Catoni, 2012]{catoni2012challenging}
Catoni, O. (2012).
\newblock Challenging the empirical mean and empirical variance: a deviation
  study.
\newblock In {\em Annales de l'IHP Probabilit{\'e}s et statistiques},
  volume~48, pages 1148--1185.

\bibitem[Chazottes et~al., 2007]{chazottes2007concentration}
Chazottes, J.~R., Collet, P., K{\"u}lske, C., and Redig, F. (2007).
\newblock Concentration inequalities for random fields via coupling.
\newblock {\em Probability Theory and Related Fields}, 137(1-2):201--225.

\bibitem[Chernoff, 1952]{chernoff1952measure}
Chernoff, H. (1952).
\newblock A measure of asymptotic efficiency for tests of a hypothesis based on
  the sum of observations.
\newblock {\em Annals of Mathematical Statistics}, 23(4):493--507.
\newblock This paper introduced Chernoff bounds and is a seminal work in the
  theory of concentration inequalities.

\bibitem[Cramér, 1938]{cramer1938mathematical}
Cramér, H. (1938).
\newblock {\em On the Mathematical Theory of Probability}.
\newblock Almqvist \& Wiksell, Stockholm.
\newblock This book is a classic in probability theory and contains the
  foundational work on large deviations and related concentration inequalities.

\bibitem[de~Finetti, 1937]{definetti1937prevision}
de~Finetti, B. (1937).
\newblock La prévision: ses lois logiques, ses sources subjectives.
\newblock {\em Annales de l'institut Henri Poincaré}, 7:1--68.

\bibitem[Diaconis and Stroock, 1989]{diaconis1989geometric}
Diaconis, P. and Stroock, D. (1989).
\newblock Geometric bounds for exponential holding times.
\newblock {\em Annals of Probability}, 17(1):361--400.

\bibitem[Douc et~al., 2011]{douc2011consistency}
Douc, R., Moulines, E., Olsson, J., and van Handel, R. (2011).
\newblock Consistency of the maximum likelihood estimator for general hidden
  {M}arkov models.
\newblock {\em Annals of Statistics}, 39(1):474--513.

\bibitem[Fan et~al., 2021]{fan2021hoeffding}
Fan, J., Jiang, B., and Sun, Q. (2021).
\newblock Hoeffding's inequality for {M}arkov chains and its applications to
  statistical learning.
\newblock {\em Journal of Machine Learning Research}, 22(139):1--35.

\bibitem[Geyer, 1992]{geyer1992practical}
Geyer, C.~J. (1992).
\newblock Practical {M}arkov chain monte carlo.
\newblock {\em Statistical Science}, 7(4):473--483.

\bibitem[Gilks et~al., 1995]{gilks1995markov}
Gilks, W.~R., Richardson, S., and Spiegelhalter, D. (1995).
\newblock {\em {M}arkov chain Monte Carlo in practice}.
\newblock CRC press.

\bibitem[Glynn and Ormoneit, 2002]{glynn2002hoeffding}
Glynn, P.~W. and Ormoneit, D. (2002).
\newblock Hoeffding's inequality for uniformly ergodic {M}arkov chains.
\newblock {\em Statistics \& probability letters}, 56(2):143--146.

\bibitem[Goodfellow et~al., 2016]{goodfellow2016deep}
Goodfellow, I., Bengio, Y., and Courville, A. (2016).
\newblock {\em Deep Learning}.
\newblock MIT Press.

\bibitem[Gorni, 1991]{gorni1991conjugation}
Gorni, G. (1991).
\newblock Conjugation and second-order properties of convex functions.
\newblock {\em Journal of Mathematical Analysis and Applications},
  158(2):293--315.

\bibitem[Hoeffding, 1963]{hoeffding1963probability}
Hoeffding, W. (1963).
\newblock Probability inequalities for sums of bounded random variables.
\newblock {\em Journal of the American statistical association},
  58(301):13--30.

\bibitem[Hoover, 1979]{hoover1979relations}
Hoover, D.~N. (1979).
\newblock {\em Relations on probability spaces and arrays of random variables}.
\newblock PhD thesis, Institute of Mathematical Statistics, Indiana University.

\bibitem[Kallenberg, 2005]{kallenberg2005probabilistic}
Kallenberg, O. (2005).
\newblock {\em Probabilistic Symmetries and Invariance Principles}.
\newblock Springer Science \& Business Media.

\bibitem[Kato, 2013]{kato2013perturbation}
Kato, T. (2013).
\newblock {\em Perturbation theory for linear operators}.
\newblock Springer.

\bibitem[Kipnis and Varadhan, 1986]{kipnis1986central}
Kipnis, C. and Varadhan, S.~S. (1986).
\newblock Central limit theorem for additive functionals of reversible {M}arkov
  processes and applications to simple exclusions.
\newblock {\em Communications in Mathematical Physics}, 104(1):1--19.

\bibitem[Kirkpatrick et~al., 1983]{kirkpatrick1983optimization}
Kirkpatrick, S., Gelatt~Jr., C.~D., and Vecchi, M.~P. (1983).
\newblock Optimization by simulated annealing.
\newblock {\em Science}, 220(4598):671--680.

\bibitem[Kontorovich and Ramanan, 2008]{kontorovich2008concentration}
Kontorovich, A. and Ramanan, K. (2008).
\newblock Concentration inequalities for dependent random variables via the
  martingale method.
\newblock {\em Annals of Probability}, 36(6):2126--2158.

\bibitem[Kontoyiannis et~al., 2005]{kontoyiannis2005relative}
Kontoyiannis, I., Lastras-Montano, L.~A., and Meyn, S.~P. (2005).
\newblock Relative entropy and exponential deviation bounds for general
  {M}arkov chains.
\newblock In {\em Proceedings of 2005 International Symposium on Information
  Theory}, pages 1563--1567. IEEE.

\bibitem[Kontoyiannis et~al., 2006]{kontoyiannis2006exponential}
Kontoyiannis, I., Lastras-Montano, L.~A., and Meyn, S.~P. (2006).
\newblock Exponential bounds and stopping rules for mcmc and general {M}arkov
  chains.
\newblock In {\em Proceedings of the 1st International Conference on
  Performance Evaluation Methodologies and Tools}. ACM.

\bibitem[Ledoux, 2001]{ledoux2001concentration}
Ledoux, M. (2001).
\newblock {\em The Concentration of Measure Phenomenon}, volume~89 of {\em
  Mathematical Surveys and Monographs}.
\newblock American Mathematical Soc., Providence, RI.
\newblock This book provides a comprehensive treatment of the concentration of
  measure phenomenon, with applications and extensions to various areas of
  probability theory.

\bibitem[Le{\'o}n and Perron, 2004]{leon2004optimal}
Le{\'o}n, C.~A. and Perron, F. (2004).
\newblock Optimal hoeffding bounds for discrete reversible {M}arkov chains.
\newblock {\em Annals of Applied Probability}, 14(2):958--970.

\bibitem[Levin et~al., 2009]{levin2009markov}
Levin, D.~A., Peres, Y., and Wilmer, E.~L. (2009).
\newblock {\em Markov Chains and Mixing Times}.
\newblock American Mathematical Society.

\bibitem[Lezaud, 1998a]{lezaud1998chernoff}
Lezaud, P. (1998a).
\newblock Chernoff-type bound for finite {M}arkov chains.
\newblock {\em Annals of Applied Probability}, 8(3):849--867.

\bibitem[Lezaud, 1998b]{lezaud1998thesis}
Lezaud, P. (1998b).
\newblock {\em Quantitative study of {M}arkov chains by perturbation of their
  kernels (in French).}
\newblock PhD thesis, University Paul Sabatier - Toulouse III.

\bibitem[Marton, 1996a]{marton1996bounding}
Marton, K. (1996a).
\newblock Bounding $\bar{d}$-distance by informational divergence: a method to
  prove measure concentration.
\newblock {\em Annals of Probability}, 24(2):857--866.

\bibitem[Marton, 1996b]{marton1996measure}
Marton, K. (1996b).
\newblock A measure concentration inequality for contracting {M}arkov chains.
\newblock {\em Geometric \& Functional Analysis}, 6(3):556--571.

\bibitem[Marton, 1998]{marton1998measure}
Marton, K. (1998).
\newblock Measure concentration for a class of random processes.
\newblock {\em Probability Theory and Related Fields}, 110(3):427--439.

\bibitem[Marton, 2003]{marton2003measure}
Marton, K. (2003).
\newblock Measure concentration and strong mixing.
\newblock {\em Studia Scientiarum Mathematicarum Hungarica}, 40(1-2):95--113.

\bibitem[Marton, 2004]{marton2004measure}
Marton, K. (2004).
\newblock Measure concentration for euclidean distance in the case of dependent
  random variables.
\newblock {\em Annals of Probability}, 32(3B):2526--2544.

\bibitem[Miasojedow, 2014]{miasojedow2014hoeffding}
Miasojedow, B. (2014).
\newblock Hoeffding’s inequalities for geometrically ergodic {M}arkov chains
  on general state space.
\newblock {\em Statistics \& Probability Letters}, 87(1):115--120.

\bibitem[Paulin, 2015]{paulin2015concentration}
Paulin, D. (2015).
\newblock Concentration inequalities for {M}arkov chains by marton couplings
  and spectral methods.
\newblock {\em Electronic Journal of Probability}, 20(79):1--32.

\bibitem[Rabiner, 1989]{rabiner1989introduction}
Rabiner, L.~R. (1989).
\newblock An introduction to hidden markov models.
\newblock {\em IEEE ASSP Magazine}, 7(2):4--16.

\bibitem[Redig and Chazottes, 2009]{redig2009concentration}
Redig, F. and Chazottes, J.~R. (2009).
\newblock Concentration inequalities for {M}arkov processes via coupling.
\newblock {\em Electronic Journal of Probability}, 14(40):1162--1180.

\bibitem[Robert et~al., 1999]{robert1999monte}
Robert, C.~P., Casella, G., and Casella, G. (1999).
\newblock {\em Monte Carlo statistical methods}, volume~2.
\newblock Springer.

\bibitem[Rockafellar, 2015]{rockafellar2015convex}
Rockafellar, R.~T. (2015).
\newblock {\em Convex analysis}.
\newblock Princeton University Press.

\bibitem[Rosenthal, 2003]{rosenthal2003asymptotic}
Rosenthal, J.~S. (2003).
\newblock Asymptotic variance and convergence rates of nearly-periodic {M}arkov
  chain monte carlo algorithms.
\newblock {\em Journal of the American Statistical Association},
  98(461):169--177.

\bibitem[Samson, 2000]{samson2000concentration}
Samson, P.-M. (2000).
\newblock Concentration of measure inequalities for {M}arkov chains and
  $\phi$-mixing processes.
\newblock {\em Annals of Probability}, 28(1):416--461.

\bibitem[Sun et~al., 2020]{sun2020adaptive}
Sun, Q., Zhou, W.-X., and Fan, J. (2020).
\newblock Adaptive {H}uber regression.
\newblock {\em Journal of the American Statistical Association},
  115(529):254--265.

\bibitem[Talagrand, 1996]{talagrand1996new}
Talagrand, M. (1996).
\newblock New concentration inequalities in product spaces.
\newblock {\em Inventiones mathematicae}, 126(3):505--563.

\bibitem[Vovk et~al., 2005]{vovk2005algorithmic}
Vovk, V., Gammerman, A., and Shafer, G. (2005).
\newblock {\em Algorithmic learning in a random world}, volume~29.
\newblock Springer.

\end{thebibliography}

\end{document}